\documentclass[12pt]{article}
\usepackage{amsmath,amsfonts,amssymb,amsthm}
\newcommand{\1}{{{\bf 1}}}
\newcommand{\id}{{\rm id}}

\newcommand{\End}{{\rm End}\,}
\newcommand{\Res}{{\rm Res}\,}
\newcommand{\Aut}{{\rm Aut}\,}

\newcommand{\haru}[2]{{\rm span}\{\,#1\,|\,#2\,\}}

\newcommand{\w}{{\omega}}

\newcommand{\Z}{\mathbb{Z}}
\newcommand{\Q}{\mathbb{Q}}
\newcommand{\C}{\mathbb{C}}

\newcommand{\h}{\mathfrak{h}}

\def\a{\alpha}
\def\b{\beta}
\def\e{\epsilon}
\def\l{\lambda}

\newcommand{\wt}[1]{{\rm wt}(#1)}

\newtheorem{theorem}{Theorem}[section]
\newtheorem{proposition}[theorem]{Proposition}
\newtheorem{lemma}[theorem]{Lemma}

\theoremstyle{definition}
\newtheorem{definition}[theorem]{Definition}

\theoremstyle{remark}
\newtheorem{remark}[theorem]{Remark}

\numberwithin{equation}{section}

\newcommand{\Free}[1]{{M(1)}^{#1}}
\newcommand{\Fremo}[1]{M(1,#1)}
\newcommand{\charge}[1]{V_{L}^{#1}}
\newcommand{\charhalf}[2]{V_{#1+L}^{#2}}
\newcommand{\charlam}[1]{V_{#1+L}}
\newcommand{\del}{\partial}

\makeatletter
\@addtoreset{equation}{section}

\makeatother

\newcommand{\Fretw}[1]{M(1)(\theta)^{#1}}
 
\newcommand{\NO}{\,{\raise0.25em
\hbox{$\mathop{\hphantom{\cdot}}\limits^{_{\circ}}_{^{\circ}}$}}\,}

\begin{document}


\begin{large}
\begin{center}
\textbf{Classification of irreducible modules \\for the vertex operator algebra $V_L^+$: General case}
\end{center}
\end{large}

\vskip 2ex
\begin{center}
Toshiyuki Abe\footnote{On leave of absence from Department of Pure and Applied Mathematics, Graduate School of Information Science and Technology, Osaka University,Toyonaka, Osaka 560-0043, Japan. Supported by JSPS Research Fellowships for Young 
Scientists.}\ \ \ and \ \ \  Chongying Dong\footnote{Partially
supported by NSF grant DMS-9987656 and a research grant from the
Committee on Research, UC Santa Cruz.}\\

\vskip 2ex  
Department of Mathematics,\\
University of California, Santa Cruz, CA 95064
\end{center}

\begin{abstract}
The irreducible modules for 
the fixed point vertex operator subalgebra $V_L^+$ 
of the vertex operator algebra
$V_L$ associated to an arbitrary positive definite even lattice $L$ 
under the automorphism lifted from the $-1$ isometry
of $L$ are classified.
\end{abstract}

\baselineskip3ex
\section{Introduction}

This paper is a continuation of our study of the $\Z_2$-orbifold
vertex operator algebra $V_L^+$ which is the fixed points of
vertex operator algebra $V_L$ (see [B], [FLM2])
associated to positive definite even lattice $L$ under the automorphism
$\theta$ lifted from the $-1$ isometry of $L.$ We classify 
the irreducible modules for $V_L^+.$ It turns out that any irreducible
$V_L^+$-module is either isomorphic to a submodule of
an irreducible 
$V_L$-module or a submodule of an irreducible twisted $V_L$-module. 
In the case that the rank of $L$ is one, this result was obtained previously 
in \cite{DN2}.

The $V_L^+$ forms a basic class of vertex operator algebras besides
affine, Virasoro and lattice vertex operator algebras.  The structure
and representation theory for affine, Virasoro and lattice vertex
operator algebras are well understood (see \cite{B}, \cite{FLM2},
\cite{D1}, \cite{DL}, \cite{DLM2}, \cite{FZ}, \cite{L1}, \cite{W})
with the help of affine Kac-Moody Lie algebras, Virasoro algebras and
lattices.  Although $V_L^+$ is also related to the lattice, but its
structure and representation theory are much more complicated. For
example, if $L$ does not have any vector of squared length 2, such as
the Leech lattice, $V_L^+$ has no weight one vectors and is hardly
related to affine Kac-Moody algebras.  In this case the weight two
subspace of $V_L^+$ is a commutative nonassociative algebra similar to
the Griess algebra \cite{Gr}.  In some sense, $V_L^+$ is a ``typical''
vertex operator algebra.  So the study of $V_L^+$ should help people
to investigate general vertex operator algebras.

It is worthy to point out that $V_L^+$ is a $\Z_2$-orbifold 
vertex operator algebra.  In the case that $L=\Lambda$ is the Leech lattice,
such vertex operator algebra was used in [FLM2] to construct the moonshine
vertex operator algebra $V^{\natural}.$  There have been efforts on 
studying the representations of $V_L^+$ for special $L$ 
with the help of representation theory for the Virasoro algebra
with central charge $1/2$ (see \cite{D3}, \cite{DGH}). We refer the
reader to [DLM1], [DM], [DVVV] and [K] for general orbifold theory.

As in \cite{DN2}, the main idea is to determine Zhu's algebra $A(V_L^+)$
whose inequivalent simple modules have a one to one correspondence
with the inequivalent irreducible (admissible) modules for $V_L^+$
(see \cite{Z}).  
The results and methods from \cite{DN2} in the case of rank one
play important roles in this paper. 
There is a fundamental difference
between rank one case and higher rank cases. Zhu's algebra in rank
one case is commutative with one exception
but is not for higher ranks. 
In general $L$ is not a orthogonal sum of rank one lattices. 
Nevertheless, 
for any element in $L$ we can consider the rank one sublattice generated
by the element and apply results in \cite{DN2}. Note that the $\theta$-invariants
$M(1)^+$ of the Heisenberg vertex operator algebra $M(1)$
is a subalgebra of $V_L^+.$ 
So the structure and representation
theory for $M(1)^+$ obtained in \cite{DN3} have been extensively used in
this paper. 
Although the information that we get on $A(V_L^+)$ is 
good enough to classify the simple modules, we could not prove that
$A(V_L^+)$ is semisimple.

There are still two major problems about $V_L^+.$ 
One is the rationality
and the other is to determine the fusion rules among irreducible
modules. 
It is certainly believed that $V_L^+$ is rational. In 
the case that $L$ is some rank one lattice, the rationality was
established in \cite{A1}. 
But rationality for any $L$ remains open. 
This problem is related to $C_2$-cofiniteness of $V_L^+,$ which
is obtained recently in [Y] in the rank one case and in \cite{ABD} in general.
The determination of the fusion rules will  be carried out in
a separate paper.

The paper is organized as follows: In Section 2, we review various
(twisted) modules for a vertex operator algebra $V$ and define Zhu's
algebra $A(V).$ 
We define vertex operator algebra $V_L^+$ and
its subalgebra $M(1)^+$ in Section 3.
 We also list the known
irreducible modules for $V_L^+.$  Section 4 is about Zhu's
algebras $A(M(1)^+)$ and $A(V_L^+).$ 
In particular we obtain
many identities for a set of generators for in $A(V_L^+)$ by using
results from \cite{DN1}, \cite{DN2} and \cite{DN3}. 
Using the known irreducible
$V_L^+$-modules we compute how the generators of $A(V_L^+)$ act
on the known simple $A(V_L^+)$-modules. 
In Section 4 we give 
a spanning set for $A(V_L^+)$ which is used in the later sections
for the classification.  
Sections 5 and 6 are devoted to the 
classification of simple $A(V_L^+)$-modules and irreducible
$V_L^+$-modules. 
Since $M(1)^+$ is a subalgebra of $V_L^+$ there
is an algebra homomorphism from $A(M(1)^+)$ to $A(V_L^+).$
This homomorphism embeds the two $d\times d$ matrix subalgebras
of $A(M(1)^+)$ (see \cite{DN3}) into $A(V_L^+)$ where $d$ is the
rank of $L.$ 
There are two cases in the classification. 
In Section 5
we deal with the case that the simple modules are annihilated by
the two matrix subalgebras. 
Section 6 handles the other case: that is, one of 
the matrix algebras does not annihilate
the simple modules. 

Throughout the paper $\Z_{\geq0}$ is the set of nonnegative integers.

\section{Vertex operator algebras and modules}\label{SVOA}

In this section we recall the definitions of various (twisted) modules
for a vertex operator algebra (cf. \cite{B}, \cite{D2}, \cite{FFR}, \cite{FLM1}, \cite{FLM2},
\cite{DLM4}, \cite{Le}, \cite{Z}). 
We also review the Zhu's algebra $A(V)$ 
associated to a vertex operator algebra $V.$

A {\em vertex operator algebra} $V$ is a $\Z$-graded vector space
$V=\bigoplus_{n\in \Z} V_n$ equipped with a linear map 
$Y:V\to(\End V)[[z,z^{-1}]],\,a\mapsto Y(a,z)=\sum_{n\in \Z}a_nz^{-n-1}$ 
for $a\in V$ such that $\dim V_n$ is finite for all integer $n$ and that $V_n=0$
for sufficiently small integer $n$ (see \cite{FLM2}).  
There are two distinguished vectors, the {\it vacuum vector} $\1\in V_0$ and the
{\it Virasoro element} $\w\in V_2$.  By definition $Y(\1,z)=\id_{V}$,
and the component operators $\{L(n)\}$ of $Y(\w,z)=\sum_{n\in\Z}L(n)z^{-n-2}$ 
gives a representation of the Virasoro algebra on $V$ with central charge $c.$  
Each homogeneous space
$V_n\,(n\geq0)$ is an eigenspace for $L(0)$ with eigenvalue $n$.

An automorphism $g$ of a vertex operator algebra $V$ is a linear
isomorphism of $V$ satisfying $g(\w)=\w$ and $gY(a,z)g^{-1}=Y(g(a),z)$
for any $a\in V$.  
We denote by $\Aut(V)$ the group of all automorphisms of $V$.  
For a subgroup $G<\Aut(V)$ the fixed
point set $V^G=\{a\in V\,|\,g(a)=a,\,g\in G\,\}$ has a
canonical vertex operator algebra structure.  

Let $g$ be an automorphism of a vertex operator algebra $V$ of order
$T$.  Then $V$ is a direct sum of eigenspaces for $g$:
\[
V=\bigoplus_{r=0}^{T-1}V^{r},\,V^{r}=\{\,a\in V\,|\,g(a)=e^{-\frac{2 \pi ir}{T}} a\,\}.
\]
\begin{definition}{\rm 
A {\em weak $g$-twisted $V$-module} $M$ is a vector space equipped
with a linear map
\begin{align*}
Y_{M}:V&\to (\End M)\{z\},\\
a&\mapsto Y_{M}(a,z)=\sum_{n\in\Q}a_nz^{-n-1},\,a_n\in \End M
\end{align*}
such that the following conditions hold for $0\leq r\leq T-1,\,a\in
V^{r},\,b\in V$ and $u\in M$:

(1) $b_mu=0$ if $m$ is sufficiently large,

(2) $Y_{M}(a,z)=\sum_{n\in\Z+\frac{r}{T}}a_nz^{-n-1}$,

(3) $Y_{M}(\1,z)=\id _{M}$,

(4) (the twisted Jacobi identity)
\begin{align*}
\begin{split}
&z_{0}^{-1}\delta\left(\frac{z_{1}-z_{2}}{z_{0}}\right)Y_{M}(a,z_{1})Y_{M}(b,z_{2})-z_{0}^{-1}\delta\left(
\frac{z_{2}-z_{1}}{-z_{0}}\right)Y_{M}(b,z_{2})Y_{M}(a,z_{1})\\
&\quad=z_{2}^{-1}\left(\frac{z_{1}-z_{0}}{z_{2}}\right)^{-\frac{r}{T}}\delta\left(\frac{z_{1}-z_{0}}{z_{2}}\right)Y_{M}(Y(a,z_{0})b,z_{2}).
\end{split}
\end{align*}
}
\end{definition}

A weak $g$-twisted $V$-module is denoted by $(M,\,Y_{M})$, or simply by $M$.  
In the case $g$ is the identity, any weak $g$-twisted $V$-module
is called a \textit{weak $V$-module}.  
A \textit{$g$-twisted weak $V$-submodule} of a $g$-twisted weak module 
$M$ is a subspace $N$ of $M$ such that $a_nN\subset N$ hold for all 
$a\in V$ and $n\in\Q$. 
If $M$ has no $g$-twisted weak $V$-submodule except $0$ and $M$, 
$M$ is called \textit{irreducible} or \textit{simple}.

Set $Y_{M}(\w,z)=\sum_{n\in\Z}L(n)z^{-n-2}$.  
Then $\{L(n)\,|\,n\in\Z\,\}$ give a representation of the Virasoro algebra
on $M$ with central charge $c$ and 
the $L(-1)$-derivative property
\begin{equation}\label{DP1} 
Y_{M}(L(-1)a,z)=\frac{d}{dz}Y(a,z)\hbox{ for all }a\in V.
\end{equation}
holds for any $a\in V$ (see \cite{DLM2}). 
  
\begin{definition}{\rm  
An \textit{admissible $g$-twisted $V$-module} $M$ is a weak $g$-twisted $V$-module 
which has a $\frac{1}{T}\Z_{\geq 0}$-gradation 
$M=\bigoplus_{n\in\frac{1}{T}\Z_{\geq 0}}M(n)$ such that 
\begin{align}\label{AD1}
a(m)M(n)\subset M(\wt{a}+n-m-1)
\end{align}
for any homogeneous $a\in V$ and $m,\,n\in\Q$.
}
\end{definition}
In the case $g$ is the identity, any admissible $g$-twisted $V$-module is called an {\em admissible $V$-module}. 
Any $g$-twisted weak $V$-submodule $N$ of a $g$-twisted admissible $V$-module is called a {\it $g$-twisted admissible $V$-submodule} if 
$N=\bigoplus_{n\in \frac{1}{T}}\Z_{\geq 0} N\cap M(n)$.

A $g$-twisted admissible $V$-module $M$ is said to be \textit{irreducible} if $M$ has no trivial admissible weak $V$-submodule. 
When a $g$-twisted admissible $V$-module $M$ is a direct sum of irreducible admissible submodules, $M$ is called \textit{completely reducible}.

\begin{definition}{\rm
A vertex operator algebra $V$ is said to be \textit{$g$-rational} if any $g$-twisted admissible $V$-module is completely reducible. 
If $V$ is $\id_{V}$-rational, then $V$ is called \textit{rational}.
}
\end{definition} 

\begin{definition}{\rm 
A \textit{$g$-twisted $V$-module} $M=\bigoplus_{\lambda\in\C} M_{\lambda}$ 
is a $\C$-graded weak $g$-twisted $V$-module with
$M_\lambda=\{u\in M|L(0)u=\lambda u\}$ such that $M_\lambda$ is finite dimensional and for fixed $\lambda\in\C$, $M_{\lambda+n/T}=0$ for sufficiently small integer $n$.
A vector $w\in M_\l$ is called a weight vector of weight $\l,$ and we write
$\l=\wt{w}.$ }
\end{definition}

In the case $g$ is the identity, a $g$-twisted $V$-module is called a \textit{$V$-module}. 
A $V$-module $M$ is called \textit{irreducible} if $M$ is irreducible as a weak $V$-module.
By definition, for a vertex operator algebra $(V,\,Y,\,\1,\,\w)$, $(V,\,Y)$ becomes a $V$-module. 
If a vertex operator algebra is irreducible as a $V$-module, then $V$ is called \textit{simple}. 

It is proved in \cite{DLM4} that if $V$ is $g$-rational then there
are only finitely many inequivalent irreducible admissible $g$-twisted
$V$-modules and any irreducible  admissible $g$-twisted
$V$-module is ordinary. 

We next define Zhu's algebra $A(V)$ which is an associative algebra following
\cite{Z}.  For any homogeneous vectors $a\in V$, and $b\in V$, we define 
\begin{align*}
&a*b=\left(\Res{z}\frac{(1+z)^{\wt{a}}}{z}
Y(a,z)\right)b,\\ 
&a\circ b=\left( \Res{z}\frac{(1+z)^{\wt{a}}}{z^{2}}
Y(a,z)\right)b,
\end{align*}
and extend to $V\times V$ bilinearly.  Denote by $O(V)$ the linear
span of $a\circ b$ ($a,b\in V$) and set $A(V)=V/O(V)$.  We write $[a]$
for the image $a+O(V)$ of $a\in V$.  The following theorem is due to
\cite{Z} (also see \cite{DLM4}).
\begin{theorem}\label{P3.1}  
(1) The bilinear operation $*$ induces $A(V)$ an associative algebra 
structure. The vector $[\1]$ is the identity and $[\w]$ is in the center of 
$A(V)$.

(2) Let $M=\bigoplus_{n=0}^{\infty}M(n)$ be an admissible $V$-module with
$M(0)\ne 0.$  
Then the linear map 
\[
o:V\rightarrow\End M(0),\;a\mapsto o(a)|M(0)
\] 
induces an algebra homomorphism from $A(V)$ to $\End M(0)$. 
Thus $M(0)$ is a left $A(V)$-module.

(3) The map $M\mapsto M(0)$ induces a bijection from the set of equivalence classes of irreducible admissible $V$-modules to the set of equivalence classes of irreducible $A(V)$-modules.   
\end{theorem}

\section{Vertex operator algebras $M(1)^+$ and $V_L^+$}

In this section we recall the construction of the vertex operator
algebra $V_L^+$ associated with a positive definite even lattice
$L$ and the vertex operator algebra $M(1)^+$ (cf. \cite{FLM2}).
We also state several results 
on classifications of irreducible modules
for $V_L^+$ when the rank of $L$ is one (see  \cite{DN2})
and $M(1)^+$ (see \cite{DN1} and \cite{DN3}).

Let $L$ be a rank $d$ even lattice with a positive definite symmetric $\Z$-bilinear form $(\cdot\,,\cdot)$.
We set $\h=\C\otimes_{\Z} L$ and extend $(\cdot\,,\cdot)$ to a $\C$-bilinear form on $\h$. 
Let $\hat{\h}=\C[t,t^{-1}]\otimes\h\oplus\C C$ be the affinization of commutative Lie algebra $\h$ defined by 
\begin{align*}
[\beta_1\otimes t^{m},\,\beta_2\otimes t^{n}]=m(\beta_1,\beta_2)\delta_{m,-n}C\hbox{ and }[C,\hat{\h}]=0
\end{align*} 
for any $\beta_i\in\h,\,m,\,n\in\Z$. 
Then $\hat{\h}^+=\C[t]\otimes\h\oplus\C C$ is a commutative subalgebra.
For any $\lambda\in\h$, we can define a one dimensional $\hat{\h}^+$-module $\C e^\lambda$ by the actions $\rho(h\otimes t^{m})e^\lambda=(\lambda,h)\delta_{m,0}e^\lambda$ and $\rho(C)e^\lambda=e^\lambda$ for $h\in\h$ and $m\geq0$.
Now we denote by 
\begin{align*}
M(1,{\lambda})=U(\hat{\h})\otimes_{U(\hat{\h}^+)}\C e^\lambda\cong S(t^{-1}\C[t^{-1}])
\end{align*}
the $\hat{\h}$-module induced from $\hat{\h}^+$-module $\C e^\lambda$.
Set $M(1)=M(1,0).$ Then there exists a linear map $Y:M(1)\to(\End
M(1,\lambda)[[z,z^{-1}]]$ such that $(M(1),\,Y,\,\1,\,\w)$ has a simple vertex
operator algebra structure and $(M(1,\lambda),Y)$ becomes an
irreducible $M(1)$-module for any $\lambda\in\h$ (see \cite{FLM2}).
The vacuum vector and the Virasoro element are given by $\1=e^0$ and
$\w=\frac{1}{2}\sum_{a=1}^{d}h_a(-1)^2\otimes e^{0}$ respectively,
where $\{h_a\}$ is an orthonormal basis of $\h$.

Let $L$ be any positive-definite even lattice and 
let $\hat{L}$ be the canonical central extension of $L$ by the cyclic
group $\langle\kappa\rangle$ of order 2:
\[
1\rightarrow \langle \kappa\rangle\rightarrow \hat{L}\mathop{\rightarrow}\limits^{-}L\rightarrow 0
\]
with the commutator map $c(\alpha,\beta)=\kappa^{(\alpha,\,\beta)}$ for
$\alpha,\beta \in L$. Let $e:L\to \hat L$ be a section such that
$e_0=1$ and $\epsilon:L\times L\to \langle\kappa\rangle$ be
the corresponding 2-cocycle. We can assume that $\e$ is bimultiplicative.
Then $\epsilon(\alpha,\beta)\epsilon(\beta,\alpha)=\kappa^{(\a,\,\b)},$
\begin{equation}\label{2c}
\e(\a,\b)\e(\a+\b,\gamma)=\e(\b,\gamma)\e(\a,\b+\gamma)
\end{equation} 
and $e_{\a}e_{\b}=\e(\a,\b)e_{\a+\b}$ for $\a,\b,\gamma\in L.$

Let $L^\circ=\{\,\lambda\in\h\,|\,(\alpha,\lambda)\in\Z\,\}$ be the
dual lattice of $L$.  
Then there is an $\hat{L}$-module
structure on $\C[L^\circ]=\bigoplus_{\lambda\in L^\circ}\C e^\lambda$ 
such that $\kappa$ act as $-1$ (see \cite{DL}). 
Let $L^{\circ}=\cup_{i\in L^{\circ}/L}(L+\lambda_i)$ be
the coset decomposition such that $\lambda_0=0.$
Set $\C[L+\lambda_i]=\bigoplus_{\alpha\in L}\C e^{\alpha+\lambda_i}.$ Then
$\C[L^\circ]=\bigoplus_{i\in L^{\circ}/L}\C[L+\lambda_i]$
and each $\C[L+\lambda_i]$ is an $\hat L$-submodule of $\C[L^\circ].$
The action of $\hat L$ on $\C[L+\lambda_i]$ is as follows:
$$e_{\alpha}e^{\beta+\lambda_i}=\e(\a,\b)e^{\a+\b+\l_i}$$
for $\alpha,\,\b\in L.$ 
On the surface, the module structure on
each $\C[L+\lambda_i]$ depends on the choice of $\lambda_i$ in
$L+\lambda_i.$ 
It is easy to prove that different choices
of $\lambda_i$ give isomorphic $\hat L$-modules.

We can identify
$e^{\alpha}$ with $e_{\alpha}$ for $\alpha\in L.$ 
Set $\C[M]=\bigoplus_{\lambda\in M}\C e^{\lambda}$ for a subset $M$ of
$L^{\circ}$, and define $V_M=M(1)\otimes\C[M]$.
Then for any $\lambda\in L^\circ$, there exists a linear map
$Y:V_L\to(\End \charlam{\lambda})[[z,z^{-1}]]$ such that
$(\charge{},\,Y,\,\1,\,\w)$ becomes a simple vertex operator algebra and
$(\charlam{\lambda},Y)$ is an irreducible $\charge{}$-module (see
[B] and \cite{FLM2}).  
The vertex operator $Y(h(-1)\1,\,z)$ and
$Y(e^\alpha,\,z)$ associated to $h(-1)\1$ and $e^\alpha$ are defined by
\begin{align*}
&Y(h(-1)\1,\,z)=h(z)=\sum_{n\in\Z}h(-n)z^{-n-1},\\
&Y(e^{\alpha},\,z)=\exp\left(\sum_{n=1}^{\infty}\frac{\alpha(-n)}{n} z^{n}\right)\exp\left(-\sum_{n=1}^{\infty}\frac{\alpha(n)}{n}z^{-n}\right)e_{\alpha}z^{\alpha},
\end{align*}
where $h(-n)\,(h\in\h,\,n\in\Z)$ is the action of $h\otimes t^n$ on $\charlam{\lambda}$, 
$e_\alpha$ is the left action of $\hat L$ on $\C[L^\circ]$ and 
$z^{\alpha}$ is the operator on $\C[L^\circ]$ defined by $z^{\alpha}e^\lambda=z^{(\alpha,\lambda)}e^\lambda$.

The vertex operator associated to the vector $v=\beta_1(-n_1)\cdots\beta_{r}(-n_r)e^\alpha$ for $\beta_i\in\h,\,n_i\geq1$ and $\alpha\in L$ is defined by 
\begin{align*}
Y(v,z)=\NO\del^{(n_1-1)}\beta_1(z)\cdots\del^{(n_r-1)}\beta_r(z)Y(e^\alpha,z)\NO,
\end{align*}
where $\del=\frac{1}{n!}(\frac{d}{dz})^{n}$ and the normal ordering $\NO\,\cdot\,\NO$ is an operation which reorders the operators so that $\beta(n)\,(\beta\in\h,n< 0$) and $e_{\alpha}$ to be placed to the left of $X(n),\,(X\in\h,n\geq 0$) and $z^{\alpha}$.
                                                                                                                                        
We note that $M(1)$ is contained in $V_L$ as a vertex operator subalgebra with same Virasoro element. 
For any $\lambda\in L^\circ$, $M(1)\otimes e^\lambda$ is isomorphic to $M(1,\lambda)$ as $M(1)$-modules. 
Thus $\charlam{\lambda}$ is isomorphic to $\bigoplus_{\alpha\in L}M(1,\lambda+\alpha)$ 
as $M(1)$-modules.

Now we define a linear isomorphism 
$\theta:V_{L+\lambda_i}\to V_{L-\lambda_i}$ for $i\in L^{\circ}/L$ by 
\begin{align*}
\theta(\beta_{1}(-n_{1})\beta_{2}(-n_{2})\cdots \beta_{k}(-n_{k})
e^{\alpha+\lambda_i})=(-1)^{k}\beta_{1}(-n_{1})\beta_{2}(-n_{2})\cdots \beta_{k}(-n_{k})
e^{-\alpha-\lambda_i}
\end{align*}
for $\beta_i\in\h,\,n_i\geq1$ and $\alpha\in L$
if $2\lambda_i\not\in L,$ and
\begin{multline*}
\theta(\beta_{1}(-n_{1})\beta_{2}(-n_{2})\cdots \beta_{k}(-n_{k})e^{\alpha+\lambda_i})\\
=(-1)^{k}c_{2\l_i}\e(\alpha,2\lambda_i)\beta_{1}(-n_{1})\beta_{2}(-n_{2})\cdots \beta_{k}(-n_{k})e^{-\alpha-\lambda_i}
\end{multline*}
if $2\lambda_i\in L$ where $c_{2\l_i}$ is a square root of
$\e(2\l_i,2\l_i).$ 
  Then $\theta$ defines a linear
isomorphism from $V_{L^{\circ}}$ to itself such that
 $$\theta Y(u,z)v=Y(\theta u,z)\theta v$$
for $u\in V_L$ and $v\in V_{L^{\circ}}.$ 
In particular, $\theta$ is an automorphism of $V_L$ which induces
an automorphism of $M(1).$ 

For any $\theta$-stable 
subspace $U$ of $V_{L\circ}$, let $U^\pm$ be the $\pm1$-eigenspace of $U$ for $\theta$. 
We have the following  proposition 
(see \cite{DM} and \cite{DLM1}):

\begin{proposition}\label{untwisted} 
(1) $\Free{\pm},\,\Fremo{\lambda}$ for $\lambda\in\h-\{0\}$ are irreducible
$M(1)^+$-modules, and $\Fremo{\lambda}\cong\Fremo{-\lambda}.$  

(2) $(V_{L+\lambda_i}+V_{L-\lambda_i})^{\pm}$ for $i\in L^{\circ}/L$
 are irreducible $\charge{+}$-modules. Moreover if $2\lambda_i\not\in L$ then
$(V_{L+\lambda_i}+V_{L-\lambda_i})^{\pm},$  $V_{L+\lambda_i},$ 
and $V_{L-\lambda_i}$ are isomorphic. 
\end{proposition}

Next we recall a construction of $\theta$-twisted modules for
$\Free{}$ and $\charge{}$ following \cite{FLM2} and \cite{D2}.  Denote
by $\h[-1]=\h\otimes t^{\frac{1}{2}}\C[t,t^{-1}]\oplus\C C$ the
twisted affinization of $\h$ defined by the commutation relation
\begin{align*}
[\beta_1\otimes t^{m},\beta_2\otimes t^{n}]=m(h_1,h_2)\delta_{m,-n}C\hbox{ and }[C,\hat{\h}]=0
\end{align*} 
for any $\beta_i\in\h,\,m,\,n\in\frac{1}{2}+\Z$. 
Then the symmetric algebra $\Fretw{}=S(t^{-\frac{1}{2}}\C[t^{-1}]\otimes\h)$ is
the unique irreducible $\hat{\h}$-module such that $C=1$ and
$\beta\otimes t^n\cdot 1=0$ if $t>0.$ 
This space, in fact, is an irreducible $\theta$-twisted $\Free{}$-module.

By abuse the notation we also use $\theta$ to denote
the automorphism of $\hat{L}$ defined by $\theta(e_{\alpha})=e_{-\alpha}$
and $\theta(\kappa)=\kappa.$  
 Set $K=\{a^{-1}\theta(a)\,|\,a\in\hat{L}\}$.
For any $\hat{L}/K$-module $T$ such that $\kappa$ acts by the scalar $-1$, we define $\charge{T}=\Fretw{}\otimes T$.
Then there exists a linear map $Y:\charge{}\to(\End \charge{T})[[z^{\frac{1}{2}},z^{-\frac{1}{2}}]]$ such that $(\charge{T},\,Y)$ becomes a $\theta$-twisted $\charge{}$-module (see \cite{FLM2}).  
The cyclic group $\langle\theta\rangle$
acts on $M(1)(\theta)$ and $\charge{T}$ by 
\begin{align*}
\theta(\beta_{1}(-n_{1})\beta_{2}(-n_{2})\cdots \beta_{k}(-n_{k}))=(-1)^{k}\beta_{1}(-n_{1})\beta_{2}(-n_{2})\cdots \beta_{k}(-n_{k})
\end{align*}
and
\begin{align*}
\theta(\beta_{1}(-n_{1})\beta_{2}(-n_{2})\cdots \beta_{k}(-n_{k})t)=(-1)^{k}\beta_{1}(-n_{1})\beta_{2}(-n_{2})\cdots \beta_{k}(-n_{k})t
\end{align*}
for $\beta_i\in\h,\,n_i\in \frac{1}{2}+\Z_{\geq0}$ and $t\in T$. 
We denote by $M(1)(\theta)^{\pm}$ and $\charge{T,\pm}$ the $\pm1$-eigenspace for $\theta$ of $M(1)(\theta)$ and $\charge{T}$ respectively.

Following \cite{FLM2}, let ${T_{\chi}}$ be the irreducible 
$\hat{L}/K$-module associated to a central character $\chi$ satisfying 
$\chi(\kappa)=-1.$ Then any irreducible $\theta$-twisted $\charge{}$-module is isomorphic to $\charge{T_{\chi}}$ for some central character $\chi$ with $\chi(\kappa)=-1$ (see \cite{D2}). 
By \cite{DLi}  we get 
\begin{proposition}\label{twisted}
(1) $\Fretw{\pm}$ are irreducible $\Free{+}$-modules.

(2) Let $\chi$ be a central character of $\hat{L}/K$ such that $\chi(\iota(\kappa))=-1$, 
and $T_{\chi}$ the irreducible $\hat{L}/K$-module with central character $\chi$. 
Then $\charge{+}$-modules $\charge{T_{\chi},\pm}$ are irreducible.
\end{proposition}

It is proved in \cite{DN1} and \cite{DN2} that any irreducible modules for $\Free{+}$ and $\charge{+}$ with rank one lattice $L$ is isomorphic to one of irreducible modules in Propositions \ref{untwisted} and \ref{twisted}:  
\begin{theorem}\label{freeclass} {\rm (\cite{DN1})} 
The set
\begin{equation}\label{IM1}
\{\,\Free{\pm},\Fretw{\pm},\Fremo{\lambda}(\cong \Fremo{-\lambda})\,|\,\lambda\in\h-\{0\}\,\}
\end{equation} 
gives all inequivalent irreducible $\Free{+}$-modules.
\end{theorem}
\begin{theorem}\label{chargeclass} {\rm (\cite{DN2})} 
Let $L=\Z\alpha$ be a rank one even lattice such that $(\alpha,\alpha)=2k$ with positive integer $k$. 
The set
\begin{equation}\label{IM2}
\{\,\charge{\pm},\charhalf{\alpha/2}{\pm},\charge{T_{i},\pm},\charlam{r\alpha/2k}\,|\,i=1,\,2,\,1\leq r\leq k-1\,\}
\end{equation}
gives all inequivalent irreducible $\charge{+}$-modules, where $T_{i}$ is an irreducible $L/2L$-module on which $e_\alpha$ acts by the scalar $(-1)^{i-1}$
for $i=1,2.$
\end{theorem}

\section{Zhu's algebras $A(\Free{+})$ and $A(V_{\Z\alpha}^+)$}

In this section we recall the structure of Zhu's algebras $A(\Free{+})$ and 
$A(\charge{+})$ following \cite{DN2} and \cite{DN3}. 

First we write down identities in Zhu's algebra $A(\Free{+})$.
Let $\{h_{a}\}$ be an orthonormal basis of $\h$. 
Set $\w_{a}=\w_{h_a}=\frac{1}{2}h_{a}(-1)^2\1$ and $J_{a}=J_{h_a}=h_{a}(-1)^4\1-2h_{a}(-3)h_{a}(-1)\1+\frac{3}{2}h_{a}(-2)^2\1$.
The vector $\w_a$ and $J_a$ generate a vertex operator algebra $\Free{+}$ associated to the one dimensional vector space $\C h_a$ (see \cite{DG}).
Next we set $S_{ab}(m,n)=h_a(-m)h_b(-n)\1$, and define $E^{u}_{ab},\,\bar{E}^{u}_{ab},\,E^{t}_{ab},\,\bar{E}^{t}_{ab}$ and $\Lambda_{ab}$ as follows (see \cite{DN3}); 
\begin{align*} 
E^u_{ab}&=5 S_{ab}(1,2)+25 S_{ab}(1,3)+36 S_{ab}(1,4)+16 S_{ab}(1,5)\,(a\neq b),\\
\bar{E}^u_{ba}&=S_{ab}(1,1)+14S_{ab}(1,2)+41S_{ab}(1,3)+44S_{ab}(1,4)+16 S_{ab}(1,5)\,(a\neq b),\\
E_{aa}^u&=E^{u}_{ab}*E^{u}_{ba},\\
E^t_{ab}&=-16(3 S_{ab}(1,2)+14S_{ab}(1,3)+19S_{ab}(1,4)+8 S_{ab}(1,5))\,(a\neq b),\\
\bar{E}^t_{ba}&=-16(5S_{ab}(1,2)+18 S_{ab}(1,3)+21 S_{ab}(1,4)+8 S_{ab}(1,5))\,(a\neq b),\\
E_{aa}^t&=E^{t}_{ab}E^{t}_{ba},\\
\Lambda_{ab}&=45 S_{ab}(1,2)+190 S_{ab}(1,3)+240 S_{ab}(1,4)+96 S_{ab}(1,5).
\end{align*}
It is proved in [DN3] that $[\bar{E}_{ab}^{u}]=[E_{ab}^u],\,[\bar{E}_{ab}^{t}]=[E_{ab}^t]$ and $[\Lambda_{ab}]=[\Lambda_{ba}]$ in $A(M(1)^+)$ for any $a,\,b$.
Thus we often use $\bar{E}_{ba}^u$ and $\bar{E}_{ba}^t$ for $[E_{ba}^u]$ and $[E_{ba}^t]$ respectively.

By \cite[Proposition 5.3.12]{DN3} we have 

\begin{proposition}\label{important} For any $a,\,b,\,c,\,d$, 
\begin{align*}
&[E^u_{ab}]*[E^u_{cd}]=\delta_{bc}[E^u_{ad}],\quad [E^t_{ab}]*[E^t*_{cd}]=\delta_{bc}[E^u_{ad}],\\
&[E^u_{ab}]*[E^t_{cd}]=[E^t_{cd}]*[E^u_{ab}]=0,\\
&[\Lambda_{ab}]*[E^u_{cd}]=[\Lambda_{ab}]*[E^t_{cd}]=[E^u_{cd}]*[\Lambda_{ab}]=[E^t_{cd}]*[\Lambda_{ab}]=0\,(a\neq b).
\end{align*}
\end{proposition}

\begin{remark}\label{remarkinverse}
The vectors $[S_{ab}(1,n)]\,(1\leq n\leq 5)$ can be expressed as linear combinations of $[E^{u}_{ab}],\,[E^{t}_{ab}]$ and $[\Lambda_{ab}]$ as follows (see \cite[Remark 5.1.2]{DN3});
\begin{align}
[S_{ab}(1,1)]&=[E_{ab}^u]+[E_{ba}^u]+[\Lambda_{ab}]+\frac{1}{2}[E_{ab}^t]+\frac{1}{2}[E_{ba}^t],\label{inverse1}\\
[S_{ab}(1,2)]&=-2[E_{ab}^u]-[\Lambda_{ab}]-\frac{3}{4}[E_{ab}^t]-\frac{1}{4}[E_{ba}^t],\label{inverse2}\\
[S_{ab}(1,3)]&=3[E_{ab}^u]+[\Lambda_{ab}]+\frac{15}{16}[E_{ab}^t]+\frac{3}{16}[E_{ba}^t],\label{inverse3}\\
[S_{ab}(1,4)]&=-4[E_{ab}^u]-[\Lambda_{ab}]-\frac{35}{32}[E_{ab}^t]-\frac{5}{32}[E_{ba}^t],\label{inverse4}\\
[S_{ab}(1,5)]&=5[E_{ab}^u]+[E_{ba}^u]+[\Lambda_{ab}]+\frac{315}{256}[E_{ab}^t]+\frac{35}{256}[E_{ba}^t].\label{inverse5}
\end{align}
\end{remark}

Let $A^u$ and $A^t$ be the linear subspace of $A(\Free{+})$ spanned by 
$E_{ab}^u$ and $E_{ab}^{t}$, respectively for $1\leq a,\,b\leq d.$  
Then we have (see \cite[Proposition 5.3.14]{DN3}):
\begin{proposition}\label{ideal} (1) $A^u$ are $A^t$ 
are two sided ideals of $A(\Free{+})$ and the quotient algebra
$A(M(1)^+)/(A^u+A^t)$ is commutative. 

(2) The natural actions of $A(M(1)^+)$
on $M(1)^{-}(0)$ and  $\Fretw{-}$ induce isomorphisms of algebras
from $A^u$ and $A^t$ to $\End\Free{-}(0)$ and $\End\Fretw{-}$ respectively. 
Under the basis $\{h_1(-1),\ldots,\,h_d(-1)\}$ (\,$\{h_1(-1/2),\ldots,\,h_d(-1/2)\}$ resp.\,) of $\Free{-}(0)$ (\,$\Fretw{-}(0)$ resp.\,),   
each $[E^{u}_{ab}]$ (\,$[E^{t}_{ab}]$ resp.\,) corresponds to the matrix element
$E_{ab}$ whose $(i,j)$-entry is 1 and zero elsewhere.
\end{proposition}

\begin{remark}\label{ideal2}
The ideals $A^u$ and $A^t$ are independent of
the choice of an orthonormal basis $\{h_{a}\}$.  
In particular, the units $I^u=\sum_{a=1}^{d}[E_{aa}^u]$ and 
$I^t=\sum_{a=1}^{d}[E_{aa}^t]$ of $A^u$ and $A^t$ are independent 
of the choice of an orthonormal basis.
\end{remark}

We next recall some relations (see \cite[Lemma 5.2.2 and Lemma 5.3.2]{DN3}) 
which will be used later.
\begin{proposition}\label{Prop-0} 
For any indices $a,\,b,\,c$,
\begin{align}
&[\w_a]*[E^{u}_{bc}]=\delta_{ab}[E^{u}_{bc}],\label{e-1}\\
&[E^{u}_{bc}]*[\w_a]=\delta_{ac}[E^{u}_{bc}],\label{e-2}\\
&[\w_a]*[E^{t}_{bc}]=\left(\frac{1}{16}+\frac{1}{2}\delta_{ab}\right)[E^{t}_{bc}],\label{e-3}\\
&[E^{t}_{bc}]*[\w_a]=\left(\frac{1}{16}+\frac{1}{2}\delta_{ac}\right)[E^{t}_{bc}],\label{e-4}\\
&[\w_a]*[\Lambda_{bc}]=[\Lambda_{bc}]*[\w_a]=0.
\end{align}
\end{proposition}

Set $H_a=H_{h_a}=J_a+\w_{a}-4\w_{a}*\w_a$ (see \cite{DN3}).
Then the following identities (see \cite[Proposition 6.13]{DN3}) hold:
\begin{proposition}\label{Prop-1} For distinct $a,\,b$ and $c$,
\begin{align} 
&\left(70 [H_a]+1188[\w_a]^2-585 [\w_a]+27\right)*[H_a]=0,\label{equation1}\\
&([\w_a]-1)*\left([\w_a]-\frac{1}{16}\right)*\left([\w_a]-\frac{9}{16}\right)*[H_a]=0,\label{equation2}\\
&-\frac{2}{9}[H_a]+\frac{2}{9}[H_b]=2[E_{aa}^u]-2[E_{bb}^u]+\frac{1}{4}[E_{aa}^t]-\frac{1}{4}[E_{bb}^t],\label{equation3}\\
\begin{split}\label{equation4}
&-\frac{4}{135}(2[\w_a]+13)*[H_a]+\frac{4}{135}(2[\w_b]+13)*[H_b]\\
&\hskip30ex=4([E_{aa}^u]-[E_{bb}^u])+\frac{15}{32}([E_{aa}^t]-[E_{bb}^t]), 
\end{split}\\
&[\w_b]*[H_a]=-\frac{2}{15}([\w_a]-1)*[H_a]+\frac{1}{15}([\w_b]-1)*[H_b],\label{equation5}\\
&[\Lambda_{ab}]^2=4[\w_a]*[\w_b]-\frac{1}{9}([H_a]+[H_b])-([E_{aa}^u]+[E_{bb}^u])-\frac{1}{4}([E_{aa}^t]+[E_{bb}^t]),\label{equation6}\\
&[\Lambda_{ab}]*[\Lambda_{bc}]=2[\w_b]*[\Lambda_{ac}].\label{equation7}
\end{align}
\end{proposition}

We now state the structure theorem on $A(\Free{+})$ (see \cite[Remark 5.1.3 and Proposition 5.3.15]{DN3}):

\begin{theorem}\label{theoremfree}  
Zhu's algebra $A(\Free{+})$ generated 
by $[w_{a}],\,[J_{a}]$ for $1\leq a\leq d,\,[\Lambda_{ab}]$ 
for $1\leq a\neq b\leq d$ and $[E^{u}_{ab}],\,[E^{t}_{ab}]$ for $1\leq a,\,b\leq d$. 
The quotient algebra $A(\Free{+})/(A^t+A^u)$ is commutative and is
generated by the images of $[w_{a}],\,[J_{a}]$ for $1\leq a\leq d$ and $[\Lambda_{ab}]$ 
for $1\leq a\neq b\leq d$.
\end{theorem}
By Theorems \ref{freeclass} and \ref{P3.1} any irreducible $A(\Free{+})$-module is a top level of an irreducible $\Free{+}$-module. 
We give the actions of the generators of $A(\Free{+})$ on the top level of all irreducible $\Free{+}$-modules.
\vskip3ex
\begin{center}
\begin{tabular}{c|c|c|c|c|c}
&$\Free{+}$&$\Free{-}$&$\Fremo{\lambda}$ $(\lambda\in\h-\{0\})$&$\Fretw{+}$&$\Fretw{-}$\\
\hline 
&$\1$&$h_c(-1)\1$&$e^{\lambda}$&$1$&$h_c(-1/2)1$\\ \hline
$\w_a$&$0$&$\delta_{ac}$&$(h_a,\lambda)^2/2$&$1/16$&$1/16+1/2\delta_{ac}$\\ 
$J_a$&$0$&$-6\delta_{ac}$&$(h_a,\lambda)^4-(h_a,\lambda)^2/2$&$3/128$&$3/128-3/8\delta_{ac}$\\
$H_a$&$0$&$-9\delta_{ac}$&$0$&$9/128$&$9/128-9/8\delta_{ac}$\\
$E_{ab}^u$&$0$&$\delta_{bc}h_a(-1)\1$&$0$&$0$&$0$\\
$E_{ab}^t$&$0$&$0$&$0$&$0$&$\delta_{bc}h_a(-1/2)\1$\\
$\Lambda_{ab}$&$0$&$0$&$(h_a,\lambda)(h_b,\lambda)$&$0$&$0$
\end{tabular}\\
\vskip2ex
Table 1. The actions of $\w,\,J,\, H,\,E^u_{ab},\,E^t_{ab}$ and $\Lambda_{ab}$ on top levels.
\end{center}
  
Let $L$ be a positive definite  even lattice. 
We next recall the results on the structure of Zhu's algebra 
$A(V_{\Z\alpha}^+)$ for any $\alpha\in L$. 
Let $\gamma\in\h$ such that $(\gamma,\gamma)\neq 0.$ 
Let $h$ be an orthonormal basis of $\h=\C\gamma$ and set $\w_\gamma=\frac{1}{2}h(-1)^2\1,\,J_{\gamma}=h(-1)^4\1-2h(-3)h(-1)\1+\frac{3}{2}h(-2)^2\1$ and $H_{\gamma}=J_\gamma+\w_{\gamma}-4\w_{\gamma}*\w_\gamma$. 
Clearly, $\w_\gamma,\,J_\gamma,\,H_\gamma$ are independent of
$\gamma$ in subspace $\C\gamma.$  
We use notations 
\begin{align*}
E^\alpha=e^\alpha+e^{-\alpha}\hbox{ and }F^\alpha=e^\alpha-e^{-\alpha}
\end{align*} for any $\alpha\in L$ respectively.
Then Zhu's algebra $A(V_{\Z\alpha}^+)$ is generated by $[\w_\alpha],\,[H_{\alpha}]$ and $[E^\alpha]$ (see \cite{DN2}).
The structure of $A(V_{\Z\alpha}^+)$ depends on $|\alpha|^2=(\a,\a)$ greatly. 

First we deal with the case $|\alpha|^2\neq2$ following \cite{DN2}. 

\begin{proposition}\label{Prop-2}
For any $\alpha\in L$ such that $|\alpha|^2=2k\neq2$, $A(V_{\Z\alpha}^+)$ is a semisimple, commutative algebra of dimension $k+7,$ and the following identities hold; 
\begin{align}
&[H_{\alpha}]*[E^{\alpha}]=\frac{18(8k-3)}{(4k-1)(4k-9)}\left([\w_{\alpha}]-\frac{k}{4}\right)\left([\w_{\alpha}]-\frac{3(k-1)}{4(8k-3)}\right)[E^{\alpha}],\label{equation8}\\
&\left([\w_{\alpha}]-\frac{k}{4}\right)\left([\w_{\alpha}]-\frac{1}{16}\right)\left([\w_{\alpha}]-\frac{9}{16}\right)[E^{\alpha}]=0.\label{equation9}
\end{align}
\end{proposition}

The actions of the generators of $A(V_{\Z\alpha}^+)$ on the
top levels of irreducible $V_{\Z\alpha}^+$-modules are given as follows:
\vskip3ex
\begin{center}
\begin{tabular}{c|c|c|c|c|c}
&$V_{\Z\a}^+$&$V_{\Z\a}^-$&$V_{\Z\a+\frac{r}{2k}\alpha}$ $(1\leq r\leq k-1)$&$V_{\Z\a+\frac{\alpha}{2}}^{+}$&$V_{\Z\a+\frac{\alpha}{2}}^{-}$\\
\hline 
&$\1$&$\alpha(-1)\1$&$e^{\frac{r}{2k}\alpha}$&$e^{\frac{\alpha}{2}}+e^{-\frac{\alpha}{2}}$&$e^{\frac{\alpha}{2}}-e^{-\frac{\alpha}{2}}$\\ \hline
$\w_\a$&$0$&$1$&$r^{2}/4k$&${k}/{4}$&${k}/{4}$\\ 
$J_\a$&$0$&$-6$&$(r^{2}/2k)^{2}-r^{2}/4k$&$k^{4}/4-k^{2}/4$&$k^{4}/4-k^{2}/4$\\
$H_\a$&$0$&$-9$&$0$&$0$&$0$\\
$E^\a$&$0$&$0$&$0$&$1$&$-1$\\ 
\end{tabular}
\vskip2ex
\begin{tabular}{c|c|c|c|c}
&$V_{\Z\a}^{T_{1},+}$&$V_{\Z\a}^{T_{1},-}$&$V_{\Z\alpha}^{T_{2},+}$&$V_{\Z\a}^{T_{2},-}$\\
\hline
&$t_1$&$\alpha(-1/2)t_1$&$t_2$&$\alpha(-1/2)t_2$\\ \hline
$\w_\a$&$1/16$&$9/16$&$1/16$&$9/16$\\ 
$J_\a$&$3/128$&$-45/128$&$3/128$&$-45/128$\\
$H_\a$&$9/128$&$-135/128$&$9/128$&$-135/128$\\
$E^\a$&$2^{-2k+1}$&$-2^{-2k+1}(4k-1)$&$-2^{-2k+1}$&$2^{-2k+1}(4k-1)$\\
\end{tabular}\\
\vskip2ex
Table 2. The actions of $\w_\a,\,J_\a,\, H_\a$ and $E^\alpha$ on top levels in the case $|\alpha|^2=2k\neq2$.
\end{center}

In the case $|\alpha|^2=2$, $V_{\Z\alpha}^+$ is isomorphic to the lattice vertex operator algebra $V_{\Z\beta}$ with $|\beta|^2=8$ (see \cite{DG}), and Zhu's algebra $A(V_{\Z\alpha}^+)$ is not a commutative but a semisimple associative algebra (see \cite{DLM3}). 
By using the representation theory of $A(V_{\Z\alpha}^+)$ we have the following proposition;
\begin{proposition}\label{Prop-alphaneq2}
For any $\alpha\in L$ with $|\alpha|^2=2$, $A(V_{\Z\alpha}^+)$ is a semisimple algebra generated by $[\w_\a],\,[J_\a],\,[E^\a]$ and the following identities hold;
\begin{align}
&[E^\alpha]*[E^\alpha]=4\e(\a,\a)[\w_\alpha],\label{equation10}\\
&[H_\alpha]*[E^\alpha]+[E^\alpha]*[H_\alpha]=-12[\w_\alpha]*\left([\w_{\alpha}]-\frac{1}{4}\right)*[E^{\alpha}],\label{equation11}\\
&([\w_{\alpha}]-1)*\left([\w_{\alpha}]-\frac{1}{4}\right)*\left([\w_{\alpha}]-\frac{1}{16}\right)*\left([\w_{\alpha}]-\frac{9}{16}\right)*[E^{\alpha}]=0.\label{equation12}
\end{align}
\end{proposition}
\begin{proof} 
The semisimplicity is due to \cite{DLM3}.  
Since $E^\alpha*E^\alpha=4\e(\a,\a)\w_\alpha$, \eqref{equation10} is clear.  
Identities \eqref{equation10}--\eqref{equation12} are immediate
as the both sides have the same actions on the simple modules
by Table 3.  
We only need to explain how to read the actions on the top
level $V_{\Z\a}^-(0)$ of $V_{\Z\a}^-$ in the table. $V_{\Z\a}^-(0)$
has a basis $\{\a(-1),\,F^\a\}$ which are eigenvectors for $[\w_a],\,[J_\a]$ and $[H_\a].$ The vector $[E^\a]$ maps $\alpha(-1)$ to $-2F^\a$ and
$F^\a$ to $2\alpha(-1).$  In order to see that $A(V_{\Z\a}^+)$ is
generated by $[E^\a],\,[\w_\a],\,[J_\a],$ we observe that the actions of
$[\w_\a],\,[J_\a],\,[E^\a]$ distinguish the simple modules.
\end{proof}

\begin{center}
\begin{tabular}{c|c|cc|c|c}
&$V_{\Z\a}^{+}$&\ \ \ \ \ \ \ $V_{\Z\a}^{-}$&&$V_{\Z\a+\frac{\alpha}{2}}^{+}$&$V_{\Z\a+\frac{\alpha}{2}}^{-}$\\
\hline 
&$\1$&$\alpha(-1)\1$&$F^\alpha$&$e^{\frac{\alpha}{2}}+c_\a e^{-\frac{\alpha}{2}}$&$e^{\frac{\alpha}{2}}-c_\a e^{-\frac{\alpha}{2}}$\\ \hline
$\w_\a$&$0$&$1$&$1$&$\frac{1}{4}$&$\frac{1}{4}$\\ 
$J_\a$&$0$&$-6$&$3$&$0$&$0$\\
$H_\a$&$0$&$-9$&$0$&$0$&$0$\\
$E^\a$&$0$&$-2F^\alpha$&$2\alpha(-1)\1$&$c_\a^3$&$-c_\a^3$\\
\end{tabular}
\vskip2ex
\begin{tabular}{c|c|c|c|c}
&$V_{\Z\a}^{T_{1},+}$&$V_{\Z\a}^{T_{1},-}$&$V_{\Z\a}^{T_{2},+}$&$V_{\Z\a}^{T_{2},-}$\\

\hline
&$t_1$&$\alpha(-1/2)t_1$&$t_2$&$\alpha(-1/2)t_2$\\ \hline
$\w_\a$&$1/16$&$9/16$&$1/16$&$9/16$\\ 
$J_\a$&$3/128$&$-45/128$&$3/128$&$-45/128$\\
$H_\a$&$9/128$&$-135/128$&$9/128$&$-135/128$\\
$E^\a$&${1}/{2}$&$-{3}/{2}$&$-{1}/{2}$&${3}/{2}$\\
\end{tabular}
\begin{center}
Table 3. The actions of $\w_\a,\,J_\a,\, H_\a$ and $E^\alpha$ on top levels in the case $|\alpha|^2=2$.
\end{center}
\end{center}

Finally we discuss the top levels of the known irreducible 
$\charge{+}$-modules. Let $\lambda\in L^{\circ}$ such that $\lambda$
has the minimal length in $L+\lambda$, i.e., $\lambda$ satisfies that $|\lambda+\alpha|^{2}\geq|\lambda|^{2}$ for any $\alpha\in L$. 
We assume that the
coset representatives $\l_i$ have the minimal lengths.
We set  
\begin{align*}
\Delta(\lambda)=\{\alpha\in L\,|\,|\,\lambda+\alpha|^2=|\lambda|^2\,\}.
\end{align*}
and $L_2=\{\,\alpha\in L\,|\,|\alpha|^2=2\,\}.$ 
We shall also use the notation
\begin{align*}
|\lambda|=\sqrt{(\lambda,\lambda)}
\end{align*}
for $\lambda\in L^{\circ}$.
Then the top levels $W(0)$ of irreducible $\charge{+}$-modules $W$ are given as follows;
\begin{align}
&\charge{+}(0)=\C\1,\quad \charge{+}(0)=\h(-1)\oplus\bigoplus_{\alpha\in L_2}\C (e^{\alpha}-e^{-\alpha}),\label{1a}\\
&\charlam{\lambda_i}(0)=\bigoplus_{\alpha\in\Delta(\lambda_i)}\C e^{\lambda_i+\alpha}\quad(2\lambda_i\notin L),\label{2a'}\\
&\charlam{\lambda_i}^{\pm}(0)=\sum_{\alpha\in\Delta(\lambda_i)}\C(e^{\lambda_i+\alpha}\pm \theta e^{\lambda_i+\alpha})\,\quad(2\lambda_i\in L),\label{3a}\\
&\charge{T_{\chi},+}(0)=T_{\chi},\quad\charge{T_{\chi},-}(0)=\h(-1/2)\otimes T_{\chi}.\label{2a}
\end{align}
Here $\h(-1)=\{h(-1)|h\in\h\}\subset M(1)$ and 
$\h(-1/2)=\{h(-1/2)|h\in\h\}\subset M(1)(\theta).$ 
\begin{remark}\label{remarklambda}
Note that if $2\lambda_i\in L$ the sum
$\sum_{\alpha\in\Delta(\lambda_i)}\C(e^{\lambda_i+\alpha}\pm \theta e^{+\lambda_i+\alpha})$ is not a direct sum since for any $\alpha\in\Delta(\lambda_i)$, 
$-2\lambda_i-\alpha$ also belongs to $\Delta(\lambda_i)$. 
Let $\bar \Delta(\lambda_i)$ be a subset
of $\Delta(\l_i)$ such that $|\bar\Delta(\lambda_i)\cap\{\a,\,-2\l_i-\a\}|=1$
for any $\alpha\in \Delta(\l_i).$
Then
\begin{align}\label{toplevelinL}
\charlam{\lambda_i}^{\pm}(0)=\bigoplus_{\alpha\in\bar\Delta(\lambda_i)}\C(e^{\lambda_i+\alpha}\pm \theta e^{\lambda_i+\alpha}).
\end{align}

\end{remark}

\section{A spanning set of $A(\charge{+})$}

For any $\alpha\in L$, set $\charge{+}[\alpha]=\Free{+}\otimes E^\alpha\oplus\Free{-}\otimes F^\alpha$ and
$A(V_L^+)(\alpha)=(V_L^+[\alpha]+O(\charge{+}))/O(\charge{+}).$ 
Then $A(V_L^+)$ is a sum of $A(V_L^+)(\alpha)$ for $\a\in L.$ 
We obtain spanning sets of $A(V_L^+)(\alpha)$ for 
$\alpha\in L$ and thus for $A(V_L^+)$ in this section. The spanning sets 
will be used in the next two sections to classify the simple 
$A(V_L^+)$-modules.   

Let $X$ be a vertex operator subalgebra of $V_L^+$ whose Virasoro
vector may differ from the Virasoro vector of $V_L^+.$ 
Clearly, the identity map induces an algebra homomorphism from
$A(X)$ to $A(V_L^+).$ 
From now on we will use $[u]$ for $u\in X$ to denote both $u+O(X)$ and $u+O(V_L^+)$ if the
context is clear. 
For example, we will use $[E_{ab}^t], [E_{ab}^u]$ for the images of $E_{ab}^t,\,E_{ab}^u$
in $A(V_L^+).$

Consider the map $\phi_\alpha:\Fremo{\alpha}\to \charge{+}[\alpha]$ defined by $u\otimes e^\alpha\mapsto\frac{1}{2}(u+\theta(u))\otimes E^\alpha+\frac{1}{2}(u-\theta(u))\otimes F^\alpha.$ Clearly
$\phi_\a$ is an $\Free{+}$-module isomorphism. 
 
First we consider a spanning set of $\Fremo{\alpha}$ for $\alpha\in \h.$ 
\begin{proposition}\label{prop-span} $\Fremo{\alpha}$ is spanned by $u(-1)e^\alpha$ and $u(-1)h(-n)e^\alpha$ for $u\in\Free{+},\,h\in\h$ and $n\geq1$. 
\end{proposition} 
\begin{proof}
Let $U$ be the subset of $\Fremo{\alpha}$  spanned by $u(-1)e^\alpha$ and $u(-1)h(-n)e^\alpha$ for $u\in\Free{+},\,h\in\h$ and $n\geq1$.
We prove that 
\begin{align}\label{inclusion}
\beta_1(-n_1)\cdots\beta_r(-n_r)e^\alpha\in U
\end{align} 
for $r\geq1,\,n_i\geq1$ and $\beta_i\in\h$ by induction on the length $r$. 
The case $r=1$ is clear.
Let $p\geq2$ and suppose that \eqref{inclusion} holds if $r<p$.  
We note that
\begin{align*}
(\gamma_1(-n_1)\cdots\gamma_r(-n_r)\1)(-1)u=&\sum\left(\prod_{j=1}^{r}\binom{-i_j-1}{n_j-1}\right)\NO \gamma_1(i_1)\cdots\gamma_r(i_r)\NO u
\end{align*} 
for any even integer $r\geq1,\,n_i\geq1,\,\gamma_i\in\h$ and $u\in\Fremo{\alpha}$, where $i_j\,(1\leq j\leq r)$ run through integers satisfying $\sum i_j=-\sum n_j$.
We see that for any $i_j\in\Z\,(1\leq j\leq r)$ such that $\sum i_j=-\sum n_j$, the coefficient $\prod_{j=1}^{r}\binom{-i_j-1}{n_j-1}$ is nonzero if and only if $i_j=-n_j$ or some $i_j$ are nonnegative.
Thus if $p$ is even then we see that 
\begin{align*}
\beta_1(-n_1)\cdots\beta_p(-n_p)e^\alpha=(\beta_1(-n_1)\cdots\beta_p(-n_p)\1
)(-1)e^\alpha+\hbox{ lower length terms},
\end{align*} 
and if $p$ is odd then  
\begin{multline*}
\beta_1(-n_1)\cdots\beta_p(-n_p)e^\alpha\\
=(\beta_1(-n_1)\cdots\beta_{p-1}(-n_{p-1})\1
)(-1)\beta_p(-n_p)e^\alpha+\hbox{ lower length terms}.
\end{multline*}
By induction hypothesis, $\beta_1(-n_1)\cdots\beta_p(-n_p)e^\alpha\in U.$
\end{proof}

By Proposition \ref{prop-span} and using the map $\phi_{\alpha},$
we see that $\charge{+}[\alpha]$ is spanned by the vectors $u(-1)E^\alpha$ and $u(-1)h(-n)F^\alpha$ for $u\in\Free{+},\,h\in\h$ and $n\geq1$.

Since 
\[
u(-1)v=u*v+\hbox{ lower weight vectors},
\]
for any homogeneous vectors $u\in\Free{+},\,v\in\charge{+}[\alpha]$ 
one proves that $\charge{+}[\alpha]$ is spanned by the vectors
$u*E^\alpha$ and $u*h(-n)F^\alpha$ for $u\in\Free{+},\,h\in\h$ and
$n\geq1$ by using induction on weight.

\begin{lemma}\label{prop-spanII}
The $A(V_L^+)(\alpha)$ is spanned by the vectors $[u]*[E^\alpha]*[v]$ for $u,\,v\in\Free{+}$.
\end{lemma} 
\begin{proof}
Let $h\in \h$ and $n>0.$ 
It is enough to prove that $h(-n)F^\a$ lies in the span
of $[u]*[E^\alpha]*[v]$ for $u,\,v\in\Free{+}$. 
If $h\in\C\a$, $h(-n)F^\a$ is contained in the spanning set by Proposition
3.15 (1) of \cite{A1}. 
Now we assume that $(h,\a)=0.$  Then 
\begin{align*}
\a(-1)h(-n)\1*E^\alpha-E^\alpha*\a(-1)h(-n)\1=(\alpha,\a)h(-n)F^\alpha.
\end{align*}
Again $h(-n)F^\a$ is contained in the spanning set. 
\end{proof}

Next we will reduce the size of the  spanning set of $V_L^+[\alpha]$ further. 
Fix an orthonormal basis $\{h_a\}$ of $\h$ so that $h_1\in\C\alpha$.
Then by Propositions \ref{Prop-2} and \ref{Prop-alphaneq2}, 
we see that $[E^\alpha]*[\w_b],\,[E^\alpha]*[H_b]$ for any $b$ are linear combinations of vectors of the form $[u]*[E^\alpha]$ for $u\in\Free{+}$. 
Since $\Free{+}$ is generated by $[\w_a],\,[H_a],\,[\Lambda_{ab}]$ and $A^u,\,A^t$,
by Lemma \ref{prop-spanII} and the fact that $[\Lambda_{ab}]$ and $[E^\alpha]$ commute if $a\neq1,\,b\neq1$,  we have:
\begin{align}\label{prespannalpha}
\begin{split}
A(V_L^+)(\a)=&\haru{[u]*[E^\alpha]}{u\in\Free{+}}\\
&+[E^\alpha]*A^u+[E^\alpha]*A^t+\sum_{b=2}^{d}\C[E^\alpha]*[\Lambda_{1b}]
\end{split}
\end{align}
where $A^u$ ($A^t$ resp.) is understood to be 
the subalgebra of $A(V_L^+)$ generated by $[E_{ab}^u]$
($[E_{ab}^t]$ resp.) for $1\leq a,b\leq d.$ Since $V_L^-(0)$ contains
$\h(-1)$ as a subspace and $V_{L}^{T_{\chi},-}$ contains
$\h(-1/2)$ as a tensor factor for any $\chi$ (see (\ref{1a}) and (\ref{2a})),
by Proposition \ref{ideal}, $A^u$ and $A^t$ are isomorphic to the matrix algebra
$M_{d}(\C).$

\begin{lemma}\label{commutelambda}
For any $b\neq1$, 
\[
[\Lambda_{1b}]*[E^\alpha]+[E^\alpha]*[\Lambda_{1b}]\in U_{1b}*[E^\alpha]+[E^\alpha]*U_{1b},
\]
where $U_{ab}$ is the subspace of $A(V_L^+)$ linearly spanned by $[E^u_{ab}],\,[E^u_{ba}],\,[E^u_{ab}]$ and $[E^u_{ba}]$.  
\end{lemma}

\begin{proof}
By Remark \ref{remarkinverse} we see that $[S_{1b}(1,1)]\equiv[\Lambda_{1b}]\equiv-[S_{1b}(1,2)]\mod U_{1b}$.
By direct calculations, we have 
\begin{multline}\label{S1n}
|\alpha|S_{1b}(1,n)*E^\alpha\\=|\alpha|h_1(-1)h_b(-n)E^\alpha+n|\alpha|^2h_b(-n-1)F^\alpha+(n+1)|\alpha|^2h_b(-n)F^\alpha.
\end{multline}
The identity 
$|\alpha|h_1(-1)h_b(-n)E^\alpha=-nh_b(-n-1)F^\alpha+L(-1)h_b(-n)F^\alpha$
shows that $|\alpha|[h_1(-1)h_b(-n)e^\alpha]=-n[h_b(-n-1)F^\alpha]-\left(n+\frac{|\alpha|^2}{2}\right)[h_b(-n)F^\alpha]$. 
Substituting this identity into \eqref{S1n} we get
\begin{multline}\label{S1nalpha-equiv}
|\alpha|[S_{1b}(1,n)]*[E^\alpha]\\=n(|\alpha|^2-1)[h_b(-n-1)F^\alpha]+\left(n(|\alpha|^2-1)+\frac{|\alpha|^2}{2}\right)[h_b(-n)F^\alpha].
\end{multline}
Similarly, 
\begin{multline}\label{alphaS1n-equiv}
|\alpha|[E^\alpha]*[S_{1b}(1,n)]\\=n(|\alpha|^2-1)[h_b(-n-1)F^\alpha]+\left(n(|\alpha|^2-1)-\frac{|\alpha|^2}{2}\right)[h_b(-n)F^\alpha].
\end{multline}
Thus, 
\begin{align}
[S_{1b}(1,1)]*[E^\alpha]-[E^\alpha]*[S_{1b}(1,1)]&=|\alpha|[h_b(-1)F^\alpha]\label{S1n-equiv},\\
[S_{1b}(1,2)]*[E^\alpha]-[E^\alpha]*[S_{1b}(1,2)]&=|\alpha|[h_b(-2)F^\alpha].\label{S2n-equiv}
\end{align}
Since $[S_{1b}(1,1)]+[S_{1b}(1,2)]\in U_{1b}$, we see that $[h_b(-2)F^\alpha]+[h_b(-1)F^\alpha]$ lies in $U_{1b}*[E^\alpha]+[E^\alpha]*U_{1b}$.
Using \eqref{S1nalpha-equiv} and \eqref{alphaS1n-equiv} gives
\begin{align*}
|\alpha|([\Lambda_{1b}]*[E^\alpha]+[E^\alpha]*[\Lambda_{1b}])&\equiv|\alpha|([S_{1b}(1,1)]*[E^\alpha]+[E^\alpha]*[S_{1b}(1,1)])\\
&\mod U_{1b}*[E^\alpha]+[E^\alpha]*U_{1b}\\
&=(|\alpha|^2-1)([h_b(-2)F^\alpha]+[h_b(-1)F^\alpha])\\
&\equiv0 \mod U_{1b}*[E^\alpha]+[E^\alpha]*U_{1b},
\end{align*} 
as desired.
\end{proof}

By Lemma \ref{commutelambda} and \eqref{prespannalpha} we immediately have
\begin{proposition}\label{spannalpha1}
For any $\alpha\in L$, 
\begin{align*}
A(V_L^+)(\a)&=\haru{[u]*[E^\alpha]}{u\in\Free{+}}
+[E^\alpha]*A^u+[E^\alpha]*A^t\\
&=\haru{[E^\alpha]*[u]}{u\in\Free{+}}
+[E^\alpha]*A^u+[E^\alpha]*A^t.
\end{align*}
In particular, any vector in $A(\charge{+})$ is a linear combination of vectors of the form $[u]*[E^\alpha]$ and $[E^\alpha]*a$ for $\a\in L,$ 
 $u\in \Free{+}$ and $a\in A^u+A^t$.
\end{proposition}

We remark that the second spanning set of $A(V_L^+)(\alpha)$ is proved 
similarly.

We conclude this section with the following lemma which will be used
in the next two sections.

\begin{lemma}\label{lemma2}
Let $I^t$ be the unit of the simple algebra $A^t$. Then for any $\alpha\in L$, $I^t*[E^\alpha]-[E^{\alpha}]*I^t=0$.
\end{lemma}
\begin{proof}

We have already pointed out that $I^t$ is independent of the
choice of orthonormal basis of $\h.$  Take an orthonormal basis $\{h_a\}$ so that $h_1\in\C\alpha$.
Then $I^t=\sum_{a=1}^{d}[E^t_{aa}]$.
It is clear that $[E^t_{aa}]*[E^\alpha]-[E^{\alpha}]*[E^t_{aa}]=0$ for any
$a\geq 2$. Hence we have to show that $[E^t_{11}]*[E^\alpha]-[E^{\alpha}]*[E^t_{11}]=0$.
Set 
\[
\bar E^t_{aa}=-16(3 S_{aa}(1,2)+18 S_{aa}(1,3)+21 S_{aa}(1,4)+8 S_{aa}(1,5)).
\]
It is easy to check by definitions that $E^t_{11}=E^t_{12}*E^t_{21}=\bar{E}^t_{11}*\bar{E}^t_{22}+u+v$, for some 
 $u\in\Free{+}_{\C h_1}$ and $v\in\Free{+}_{\C h_2}$ where
$M(1)^+_{W}$ is the subspace of $M(1)^+$ corresponding to
the subspace $W$ of $\h.$ 
Since $[\bar{E}^t_{22}]$ and $[v]$ commute with $[E^\a]$ we only need
to prove that $[\bar{E}^t_{11}]$ and $[u]$ commute with $[E^\a].$ 
Note that 
 the identity map induces an algebra homomorphism from $A(V_{\Z\a}^+)$ 
to $A(V_L^+).$ It suffices to prove that  $[\bar{E}^t_{11}]$ and $[u]$ commute
with $[E^\a]$ in $A(V_{\Z\a}^+).$

The result is clear if $|\alpha|^2\neq 2$ as $A(V_{\Z\alpha}^{+})$ 
is commutative (see Proposition \ref{Prop-2}).  

If $|\alpha|^2=2,$ $A(V_{\Z\a}^+)$ is a semisimple algebra of
dimension 11 by Proposition \ref{Prop-alphaneq2} and Table 3. 
So it is enough to prove that $[\bar{E}^t_{11}]$ and $[u]$ commute with
$[E^\a]$ 
on the unique 2 dimensional module $V_{\Z\a}^-(0)$ spanned by
$\a(-1)$ and $F^\a.$ 
Since $[E_{11}^t]=0$ on $V_L^-(0)$ we immediately
see that $[E_{11}^t]$ acts on $V_{\Z\a}^-(0)\subset V_L^-(0)$ as zero.  
Clearly, $[\bar{E}^t_{22}]$ and $[v]$ act on $V_{\Z\a}^-(0)$  trivially. 
As a result, $[u]=0$ on $V_{\Z\a}^-(0).$  
A direct calculation shows that $o(\bar E^t_{11})h_1(-1)=o(\bar E^t_{11})(F^\alpha)=0$. This shows that $[\bar E^t_{11}]$ and $[u]$ are
zero on $V_{\Z\a}^-(0), $ and commute with $[E^\a]$ in particular.
So the proof is complete.
\end{proof}

\section{Classification I}

In this section and the next section we classify the simple modules
for $A(V_L^+)$ and thus classify the irreducible admissible 
$V_L^+$-modules. We prove that any irreducible admissible 
$V_L^+$-module is ordinary and is isomorphic to one given in
Propositions   \ref{untwisted} and  \ref{twisted}. In this section
we deal with a simple $A(V_L^+)$-module $W$ such that $A^uW=A^tW=0.$
The other cases will be studied in the next section.

Now let $W$ be an irreducible $A(\charge{+})$-module such that $A^uW=A^tW=0$.
An element $u\in A(V_L^+)$ is called {\em semisimple} on $W$ if $u$ acts on 
$W$ diagonally. Note that for any $\alpha\in L,$ $A(V_{\Z\a}^{+})$ 
is semisimple, and that $[\w_\a]$ and $[J_\a]$ are semisimple on $W.$ 
As a result, $[H_\a]$ is semisimple  on $W$ 
since $[\w_a]$ and $[J_a]$ commute.

Take mutually orthogonal elements $\alpha_a\in L$ 
with  $1\leq a\leq d$, and consider the orthonormal basis $\{h_a\}$ of $\h$ such that $h_a\in\C\alpha_a$. 
Since $A^tW=A^uW=0$ we see that  $[S_{ab}(1,1)]=[\Lambda_{ab}]$ on
$W.$ Using identity 
$$\frac{1}{2}(\alpha_a+\alpha_b)(-1)^{2}\1=\w_a+|\alpha_a||\alpha_b|
S_{ab}(1,1)+\w_b,$$
we have  
\begin{align}\label{lambdaabequation}
\frac{1}{2}[(\alpha_a+\alpha_b)(-1)^{2}\1]=[\w_a]+|\alpha_a||\alpha_b|[\Lambda_{ab}]+[\w_b]
\end{align}
on $W$.
Since $[\w_a],\,[\w_b]$ and $[(\alpha_a+\alpha_b)(-1)^{2}\1]$ are semisimple and commute each other on $W,\,[\Lambda_{ab}]$ is also semisimple on $W.$ 
Note that $A(M(1)^+)$ is generated by $[\w_a],\,[J_a]$ for $1\leq a\leq d$, $[\Lambda_{ab}]$ for $1\leq a\neq b\leq d$ and $A^u,\,A^t$.  
We obtain the following lemma.
\begin{lemma}\label{diagonal} 
Let $W$ be an irreducible $A(\charge{+})$-module such that $A^uW=A^tW=0$. 
Then any element in $A(\Free{+})$ is semisimple on $W.$ 
\end{lemma}

Since $A(M(1)^+)/(A^u+A^t)$ is a commutative algebra by Proposition \ref{ideal},
Lemma \ref{diagonal} implies that $W$ is a direct sum of one-dimensional
irreducible $A(\Free{+})$-modules on which $A^u$ and $A^t$ act as zero.  
By Theorem \ref{freeclass}, an irreducible $A(\Free{+})$-submodule of $W$ is
isomorphic to one of the following: $\Free{+}(0),\,\Fretw{+}(0)$ or
$\Fremo{\lambda}(0)$ for some $\lambda\in\h-\{0\}$.

Now we consider the case $[H_\gamma]W=0$ for some $\gamma\in\h$.
Then we take an orthonormal basis $\{h_a\}$ of $\h$ so that $h_1\in\C\gamma$.
By using \eqref{equation3}, we have 
\begin{align}\label{ha-equation}
[H_a]=[H_1]=0
\end{align}
on $W$ for any $a$.  
We note from Table 1 that any irreducible
$A(\Free{+})$-module on which $[H_a]$ acts $0$ for $a=1,\ldots,d$ is
isomorphic to $\Free{+}(0)$ or $\Fremo{\lambda}(0)$ for
$\lambda\in\h-\{0\}$.  
So $W$ is a direct sum of $\Free{+}(0)$ and $\Fremo{\lambda}(0)$ 
as an $A(\Free{+})$-module.  
In fact, $[H_\gamma]W=0$ for any $\gamma\in\h$ with
$(\gamma,\gamma)\ne 0$ again by Table 1.  

We write $W=W_{0}\oplus\left(\sum_{\lambda\neq 0}W_{\lambda}\right)$, 
where $W_{\lambda}$ is a sum of simple $A(\Free{+})$-submodules isomorphic
to $\Fremo{\lambda}(0)$ if $\l\ne 0$ and $W_0$ is
a sum of simple $A(\Free{+})$-submodules isomorphic
to $M(1)^+(0)=\C\1.$   
Since $\Fremo{\lambda}\cong\Fremo{-\lambda}$ as $\Free{+}$-modules, 
we can assume that $\lambda$ ranges in $\h-\{0\}/\sim$, 
where the equivalence $\lambda\sim \mu$ is defined
by $\lambda=\pm\mu$.

We first assume that $W_0\neq0$, 
and let $v$ be a nonzero vector in $W_{0}$. 
Then $\C v$ is isomorphic to $\C\1$ as $A(\Free{+})$-modules. 
This shows that $[\w_{\alpha}]v=0$ for any $\alpha\in L$. 
Hence by \eqref{equation9} and \eqref{equation12} we have $[E^{\alpha}]v=0$. 
Since $A(\charge{+})$ is generated by $[u]$ and $[E^{\alpha}]\,(u\in\Free{+},\,\alpha\in L)$, $\C v$ is an $A(\charge{+})$-module isomorphic to $\charge{+}(0)=\C\1$. 
This shows that $W$ is isomorphic to $\charge{+}(0)$.

Next we consider the case $W_{0}=0$.  
We set $P(W)=\{\lambda\in \h-\{0\}\,|\,W_{\lambda}\neq0\,\}$.

\begin{lemma}\label{prop-lambda2} (1) For any $\lambda,\,\mu\in P(W)$, $|\lambda|^2=|\mu|^2$. 

(2) For any $\alpha\in L$ and $\lambda\in P(W)$, $(\lambda,\alpha)\in\Z$, i.e., $P(W)\subset L^{\circ}$. 

(3) For any $\l\in P(W),$ $\l$ has the minimal length in $L+\l.$  
\end{lemma} 
\begin{proof}
(1) follows from the fact that the Virasoro element $[\w]$ acts on $W$ as a constant and acts on $W_{\lambda}$ by scalar $\frac{|\lambda|^2}{2}$. 

For (2) we note that for any $\alpha\in L$ and that $[\w_{\alpha}]$ acts on $W_{\lambda}$ by the scalar $\frac{(\lambda,\alpha)^2}{2|\alpha|^2}$. 
On the other hand, by Tables 2 and 3, $[\w_\a]$ acts
on any irreducible $A(V_{\Z\alpha}^+)$-module on which $[H_{\alpha}]=0$ 
as $\frac{r^2}{2|\alpha|^2}\,(0\leq r\leq \frac{|\alpha|^2}{2}).$ 
Thus for any $\lambda\in P(W)$, $\frac{(\lambda,\alpha)^2}{2|\alpha|^2}=\frac{r^2}{2|\alpha|^2}$ for some $r$. 
This shows  that $(\lambda,\alpha)$ is an integer.

We now prove (3) and let $\alpha,r$ be as in the proof of (2). 
Then $(\l,\a)^2=r^2\leq\frac{|\a|^4}{4}.$ Thus $2(\l,\a)\geq -|\a|^2$
which is equivalent to that $|\lambda+\alpha|^2\geq |\l|^2.$
\end{proof}   

For $\lambda\in P(W)$, we take a nonzero vector $v\in W_\lambda$ and set $Q(\lambda,W)=\{\alpha\in L\,|\,[E^\alpha]v\neq0\}$. 
We see that 
\[
W=\sum_{\alpha\in\Q(\lambda,W)}\C([E^\alpha]v)
\] 
by Lemma \ref{diagonal} and Proposition \ref{spannalpha1}.
By using Proposition \ref{Prop-2}, Proposition \ref{Prop-alphaneq2} and the fact that $[H_{\alpha}]=0$ on $W$, we get $\frac{|\alpha|^2}{8}=\frac{(\lambda,\alpha)^2}{2|\alpha|^2}$.
This shows that $\alpha$ satisfies 
\begin{align}\label{lambda+alpha}
|\lambda+\alpha|^2=|\lambda|^2\hbox{ or }|\lambda-\alpha|^2=|\lambda|^2.
\end{align}
That is, $\alpha\in \Delta(\lambda)\cup\Delta(-\lambda)$.

\begin{lemma}\label{qlambdaset}
For any nonzero $\lambda\in P(W)$, $Q(\lambda,W)=\Delta(\lambda)\cup\Delta(-\lambda)$.
\end{lemma}
\begin{proof}
We have already known that 
$Q(\lambda,W)\subset\Delta(\lambda)\cup\Delta(-\lambda)$. 
To prove $\Delta(\lambda)\cup\Delta(-\lambda)\subset Q(\lambda,W)$ it suffices to show that $[E^\alpha]v\neq0$ for any nonzero $\alpha\in \Delta(\lambda)$. 
Write $\lambda=\lambda_1+\lambda_2$ such that
$\lambda_1\in\C\alpha$ and $(\lambda_2,\alpha)=0$.
Then the condition $|\lambda+\alpha|^2=|\lambda|^2$ implies that 
\[
2(\lambda_1,\alpha)=-|\alpha|^2.
\]
Hence $\lambda_1=\frac{(\lambda_1,\alpha)}{|\alpha|^2}\alpha=-\frac{1}{2}\alpha$. 
By Theorem \ref{chargeclass} and Table 2 and Table 3, 
$\C v$ is an irreducible $A(V_{\Z\alpha}^+)$-module
isomorphic to $V_{\frac{\alpha}{2}+\Z\alpha}^+(0)$ or $V_{\frac{\alpha}{2}+\Z\alpha}^-(0)$.  Thus Table 2 and Table 3 again implies $[E^\alpha]v\neq0$,
as desired.  
\end{proof}  

Since $\Delta(\lambda)\cap\Delta(-\lambda)=\{0\}$ for any $\lambda\in L^\circ$
and $E^\alpha =E^{-\alpha}$ we have
\begin{align}\label{decomp-lambda}
W=\sum_{\alpha\in\Delta(\lambda)}\C([E^\alpha]v).
\end{align}
Recall from Section 4 that $\bar\Delta(\lambda)$ is a subset 
of $\Delta(\lambda)$ if $2\lambda\in L.$ 

\begin{lemma}\label{lemma-2lambda}
(1) If $2\lambda\notin L$, then $W_{\l+\alpha}=\C([E^\alpha]v)$ for
$\alpha\in \Delta(\lambda),$ $W_{\l+\alpha}=0$ for other $\alpha$ and 
$W=\bigoplus_{\alpha\in\Delta(\lambda)}\C([E^\alpha]v)$.

(2) If $2\lambda\in L$, then 
$$W_{\l+\alpha}=\C([E^\alpha]v)=\C([E^{\alpha+2\lambda}]v)$$ for
$\alpha\in \Delta(\lambda),$ $W_{\l+\alpha}=0$ for other $\alpha$ and 
$W=\bigoplus_{\alpha\in\bar\Delta(\lambda)}\C([E^\alpha]v)$. 
\end{lemma}
\begin{proof}
Let $\alpha\in \Delta(\lambda)$ and $h\in\h$. 
Then $h=\frac{(h,\alpha)}{|\alpha|^2}\alpha+(h-\frac{(h,\alpha)}{|\alpha|^2}\alpha)$.  
We take an orthogonal basis $\{h_a\}$ so that $h_1\in\C\alpha.$  
By Lemma \ref{commutelambda} we have 
\begin{align*}
[\Lambda_{1i}][E^\alpha]v=-[E^\alpha][\Lambda_{1i}]v=-(h_1,\lambda)(h_i,\lambda)[E^\alpha]v
\end{align*}
for $2\leq i\leq d$
and
\begin{align*}
[\Lambda_{ij}][E^\alpha]v=[E^\alpha][\Lambda_{ij}]v=-(h_i,\lambda)(h_j,\lambda)[E^\alpha]v
\end{align*}
for $2\leq i<j\leq d.$ 
Note that $[\w_i][E^\alpha]v=\frac{(h_i,\lambda)^{2}}{2}[E^\alpha]v$ for
all $i.$ We have
\begin{align*}
\frac{1}{2}[h(-1)^2\1][E^\alpha]v=&\sum_{i=1}^{d}(h,h_{i})^{2}[\w_{i}][E^\alpha]v+\sum_{1\leq i<j\leq d}(h,h_{i})(h,h_{j})[\Lambda_{ij}][E^\alpha]v\\
=&\sum_{i=1}^{d}\frac{(\lambda,h_{i})^{2}(h,h_{i})^{2}}{2}[E^\alpha]v-\sum_{2\leq j\leq d}(\lambda,h_{1})(\lambda,h_{j})(h,h_{1})(h,h_{j})[E^\alpha]v\\
&\quad+\sum_{2\leq i<j\leq d}(\lambda,h_{i})(\lambda,h_{j})(h,h_{i})(h,h_{j})[E^\alpha]v\\
=&\frac{1}{2}\left(((h,h_1)h_1,\lambda)-\sum_{i=2}^{d}(h,h_{i})(\lambda,h_{i})\right)^2[E^\alpha]v\\
=&\frac{1}{2}\left(2((h,h_1)h_1,\lambda)-(h,\lambda)\right)^2[E^\alpha]v\\
=&\frac{1}{2}\left(2(\frac{(h,\alpha)}{|\alpha|^2}\alpha,\lambda)-(h,\lambda)\right)^2[E^\alpha]v\\
=&\frac{1}{2}\left(-(h,\alpha)-(h,\lambda)\right)^2[E^\alpha]v\\
=&\frac{1}{2}\left(h,\lambda+\alpha\right)^2[E^\alpha]v.
\end{align*} 
This implies that $[E^\alpha]v\in U_{\lambda+\alpha}$. 
Notice that 
 $\Fremo{\lambda+\alpha}\cong \Fremo{\lambda+\beta}$ as $M(1)^+$-modules
if and only if $\alpha=\beta$ or $2\lambda=-\alpha-\beta$ for any $\alpha,\,\beta\in L$. 
Thus if $2\lambda\notin L$ then $\sum_{\alpha\in\Delta(\lambda)}U_{\lambda+\alpha}=\bigoplus_{\alpha\in\Delta(\lambda)}U_{\lambda+\alpha}$.
This proves (1). 

Next we consider the case $2\lambda\in L$. Clearly, $W=\bigoplus_{\alpha\in \bar\Delta(\lambda)}W_{\l+\a}$ and 
$W_{\l+\alpha}=\C([E^\alpha]v)+\C([E^{\alpha+2\lambda}]v)$ for 
$\alpha\in \Delta(\lambda).$ Since $\lambda\in P(W)$ is arbitrary,
it suffices to prove that $W_{\l}=\C v=\C([E^{2\lambda}]v).$ Note that
$W_{\l}$ is an $A(V_{\Z2\lambda}^+)$-module which is a sum
of $V_{\Z2\l+\l}^{\pm}(0)$ (see Tables 2 and 3). Thus we can choose
$v$ to be an eigenvector for $[E^{2\lambda}].$ As a result
we see that  $\C([E^{2\lambda}]v)=\C v,$ as desired.
\end{proof}

We now can prove the following:
\begin{proposition}\label{propositonlambda} 
Let $W$ and $\lambda$ be as before. Then
$W$ is isomorphic to $V_{L+\l}(0)$ if $2\l\not\in L$ and
to $V_{L+\l}^{\pm}$ if $2\l\in L$ and $[E^{2\l}]=\pm c^3_{2\l}$ on $W_\l.$
\end{proposition} 
\begin{proof} We first assume that $2\l\notin L.$ 
We notice that $V_{L+\l}(0)=\bigoplus_{\a\in\Delta(\l)}\C
e^{\lambda+\a}$ by Lemma \ref{prop-lambda2} (3) and \eqref{2a'} and
that $W=\bigoplus_{\a\in\Delta(\l)}\C[E^\a]v$ for a nonzero $v\in W_\l$ 
by Lemma \ref{lemma-2lambda}.  Now we define a linear isomorphism $f$ from $W$
to $V_{L+\l}(0)$ by sending $ [E^\a]v$ to $o(E^{\a})e^{\l}$ for $\a\in
\Delta(\l)$.  Since $[E^{\a}]v=o(E^{\a})e^{\l}=0$ for $\a\notin
\Delta(\l)\cup\Delta(-\l),$ we see immediately that $f$ sends $
[E^\a]v$ to $o(E^{\a})e^{\l}$ for any $\a\in L.$ Clearly, both $ [E^\a]v$
and $o(E^{\a})e^{\l}$ are eigenvectors for $A(M(1)^+)$ with the same
eigenvalue. Thus we have
$$f([u*E^{\a}]v)=f([u][E^{\a}]v)=o(u)o(E^\a)e^{\l}=o(u*E^\a)e^{\l}$$
for any $u\in M(1)^+$ and $\a\in L.$ 
Recall that both $A^u$ and $A^t$
act trivially on $W$ and $V_{L+\l}(0).$
By Proposition \ref{spannalpha1}
we see that 
$$f([u]v)=o(u)f(v)=o(u)e^{\l}$$
for any $u\in V_L^+.$ 
Let $w=[x]v$ for some $x\in V_L^+$ and
$u\in V_L^+.$ 
Then
$$f([u]w)=f([u*x]v)=o(u*x)f(v)=o(u)o(x)f(v)=o(u)f([x]v)=o(u)f(w).$$
That is, $f$ is an $A(V_L^+)$-module isomorphism.

The case that $2\l\in L$ is more complicated. In this case
$[E^{2\l}]v=\pm c_{2\l}^3v.$ We first assume that  $[E^{2\l}]v=c_{2\l}^3v.$
Note that $V_{L+\l}(0)=\bigoplus_{\a\in\bar{\Delta}(\l)}\C(e^{\lambda+\a}+\theta e^{\lambda+\a})$ with $o(E^{2\l})(e^{\lambda}+\theta e^{\lambda})=c_{2\l}^3(e^{\lambda}+\theta e^{\lambda})$ by Lemma \ref{prop-lambda2} (3) and \eqref{3a}, where we take $\lambda$ to be a representative of $\lambda+L$.
As before we define a linear map $f:W\to \charlam{\lambda}^+(0)$ by $f([E^\alpha]v)=o(E^\alpha)(e^\lambda+e^{-\lambda})$ for any 
$\alpha\in\bar\Delta(\lambda).$ 
It is clear that $f$ is an $A(\Free{+})$-module isomorphism. Using the
proof for the case that $2\l\notin L$ it is enough to
prove that $f([E^\alpha]v)=o(E^\alpha)f(v)$ for any $\alpha\in L$.
Since $o(E^\alpha)f(v)=f([E^\alpha]v)=0$ for any 
$\alpha\in L-\Delta(\lambda)\cap \Delta(-\l)$ and $E^{\a}=E^{-\a}$ 
we only need to show that $f([E^{2\lambda+\alpha}]v)=o(E^{2\lambda+\alpha})f(v)$ for any $\alpha\in\bar\Delta(\lambda).$

Let $0\neq\alpha\in \Delta(\l).$ Then $(\l,\a)$ is a negative integer.
Thus $e^{\a}_{n}e^{-2\l}=e^{-\a}_ne^{2\l}=0$ for all $n\geq -1$ and
$E^\a*E^{2\l}\in V_L^+[\a+2\l].$ 
As a result we see that
$[E^\a]*[E^{2\l}]=[u]*[E^{\a+2\l}]$ on both $W$ and
$V_{L+\l}^+(0)$ for some $u\in M(1)^+$ by 
Proposition \ref{spannalpha1}.  Since $\C
o(E^{2\lambda+\alpha})(e^\lambda+\theta e^{\lambda})\cong\C[E^{2\lambda+\alpha}]v\cong\Fremo{\lambda+\alpha}(0)$
as $A(\Free{+})$-modules, $[u]$ acts on $\C
o(E^{2\lambda+\alpha})(e^\lambda+\theta e^{\lambda})$ and
$\C[E^{2\lambda+\alpha}]v$ as a same constant $p.$ Then
$[E^\alpha]v=[E^\alpha][E^{2\lambda}]v=p[E^{2\lambda+\alpha}]v$ and
$o(E^\alpha)(e^{\l}+\theta e^{\l})= po(E^{\alpha+2\l})(e^{\l}+\theta e^{\l}).$
Since $[E^\alpha]v$ is nonzero, $p$ is nonzero. 
This implies that
$f([E^{2\lambda+\alpha}]v)=o(E^{2\lambda+\alpha})f(v).$ 
So if $[E^{2\l}]v=c_{2\l}^3v,$ $W$ is isomorphic to $V_{L+\l}^+.$

Similarly, if $[E^{2\l}]v=-c_{2\l}^3v,$ $W$ is isomorphic to $V_{L+\l}^-.$ The proof
is complete.
\end{proof} 

Next we consider the case $[H_\gamma]W\neq0$ for some nonzero $\gamma\in\h$. 
By Lemma \ref{diagonal}, there exists an eigenvector $v\in W$ for $[H_\gamma]$.
We take an orthonormal basis $\{h_a\}$ so that $h_1\in\C\gamma$.
Then \eqref{ha-equation} implies $[H_a]v=[H_b]v\neq0$ for any $a,\,b$. 
By \eqref{equation4} we have $[\w_a]v=[\w_b]v$ for any $a$ and $b$. 
Thus \eqref{equation5} implies $[w_a]v=\frac{1}{16}v$ for any $a$, and \eqref{equation1} shows that $[H_a]v=\frac{9}{128}v$. 
We finally have $[\Lambda_{ab}]v=0$ by \eqref{equation6}. 
By Theorem \ref{freeclass} and Table 1, we see that $\C v$ is isomorphic to $\Fretw{+}(0)$ as an $A(\Free{+})$-module. This shows that $H_{\beta}v=\frac{9}{128}v$ for any $\beta\in \h$ with $(\beta,\beta)\ne 0.$ Thus
$W=W_{9/128}\oplus W_0$ where $W_{9/128}$ and $W_0$ are subspaces
of $W$ such that $H_\b$ acts as $\frac{9}{128}$ and
$0$ respectively for any $\beta\in \h$ with $(\beta,\beta)\ne 0.$ 

Let $\alpha\in L.$ 
We claim that $W_{9/128}$ is $[E^{\a}]$-invariant.
It is enough to show that $[H_{\a}][E^{\a}]w\ne 0$
if $[E^\a] w$ is nonzero for $w\in W.$ 
By Propositions \ref{Prop-2} and 
\ref{Prop-alphaneq2}, $0\ne [E^{\a}]w$ is an eigenvector
for $[H_{\a}]$ with nonzero eigenvalue by noting
that $[\w_{\a}]*[E^{\a}]=[E^\a]*[\w_{\a}]$ and $[\w_\a]$ acts
on $W_{9/128}$ as constant $1/16.$  By Proposition \ref{spannalpha1}
$W_{9/128}$ is a submodule and  must be $W$ itself. 
Thus
$W$ is a direct sum of copies of $\Fretw{+}(0)$ and each element
of $A(M(1)^+)$ acts on $W$ as a constant. 
This implies that the action of $[E^\alpha]$ for any $\alpha\in L$ commutes with the action of $A(\Free{+})$.

For any $0\ne \alpha\in L$ we define 
\begin{align}\label{groupaction}
B_{\alpha}=2^{|\alpha|^2-1}\left(E^\alpha-\frac{2|\alpha|^2}{2|\alpha|^2-1}E^t_{11}*E^\alpha\right)
\end{align}
and $B_0=\1$ where $E^t_{11}$ is defined with respect to
an orthonormal basis $\{h_a\}$ of $\h$ such that $h_1\in\C\a.$
Note that $A^tW=0.$ 
So in our case, $[B_{\alpha}]=2^{|\alpha|^2-1}[E^\alpha].$
Let $\alpha,\beta\in L$ such that $(\alpha,\beta)<0$. 
The introduction
of $B_{\a}$ is motivated by the definition of vertex operator
$Y(E^{\a},z)$ on the twisted modules $V_{L}^{T_{\chi}}.$ Since
$E^\a*E^\b\in V_L^+[\a+\b],$  
by Proposition \ref{spannalpha1}, we see that  
\[
[B_{\alpha}]*[B_{\beta}]\equiv[u]*[B_{\alpha+\beta}]\mod [E^\alpha]*(A^u+A^t),
\]
for some $u\in \Free{+}$ on $W$. 
Since $[u]$ acts on $\Fretw{+}(0)$ by a constant, say $p$, we have $[B_{\alpha}]*[B_{\beta}]=p[B_{\alpha+\beta}]$ on $W$.  

On the other hand,   
\[
o(B_{\alpha}*B_{\beta})=o(u)o([B_{\alpha+\beta}])=\epsilon(\alpha,\beta)o(B_{\alpha+\beta}),
\]
on the irreducible $A(\charge{+})$-module $\charge{T_{\chi},+}(0)$ for any $\chi$. As a result we see that $p=\epsilon(\alpha,\beta)$. 

In the case $(\alpha,\beta)\geq 0$, we have
\begin{align}\label{zerocase}
\begin{split}
[B_{\alpha}]*[B_{\beta}]&=\epsilon(\alpha+\beta,\beta)[B_{\alpha+\beta}]*[B_{-\beta}]*[E^{\beta}]\\
&=\epsilon(\alpha+\beta,\beta)\epsilon(\beta,\beta)[B_{\alpha+\beta}]*[B_{0}]\\
&=\epsilon(\alpha,\beta)[B_{\alpha+\beta}].
\end{split}
\end{align}
So in any case we have $[B_{\alpha}]*[B_{\beta}]=\epsilon(\alpha,\beta)[B_{\alpha+\beta}]$ for $\a,\b\in L.$  

Thus the map from $\hat L$ to $GL(W)$ by sending
$e_{\alpha}$ to $B_{\a}$ and $\kappa$ to $-1$ defines representation
of $\hat L$ on $W.$ Note that $\theta(e_\a)=e_{-\a}$ and
$B_{\a}=B_{-\a}.$ So this induces a group homomorphism   
from $\hat{L}/K$ to $GL(W).$ Since the action
of $A(M(1)^+)$ commutes with the action of $[B_\a]$ for any $\a\in L,$  
by Proposition \ref{spannalpha1} we see that $W$ is an irreducible $\hat{L}/K$-module. Thus $W$ is isomorphic 
to $T_{\chi}$ for some central character $\chi$ with $\chi(\kappa)=-1$.
We finally obtain the following proposition: 
\begin{proposition}\label{twistplus}
Let $W$ be an irreducible $A(\charge{+})$-module such that $A^uW=A^tW=0$. 
If $[H_\gamma] W\neq 0$ for some $\gamma\in\h$, then  
there exists an irreducible $\hat{L}/K$-module $T_{\chi}$ with central character $\chi$ such that $W\cong T_{\chi}=\charge{T_{\chi},+}$. 
\end{proposition}

\section{Classification II} 

Again we fix an irreducible $A(V_L^+)$-module $W.$ 
In this section we consider the case $A^uW\neq0$ or $A^tW\neq0$ and complete
our classification of irreducible $A(V_L^+)$-modules. 

First we assume that $A^uW$ is nonzero. 
Since $A^u$ is a simple 
algebra, $W$ contains a simple 
$A^u$-modules isomorphic to $\h(-1).$ 
So we can assume that
$\h(-1)$ is an $A(M(1)^+)$-submodule of $W.$ We will use the action of
$A(V_L^+)$ on $\h(-1)$ to get the whole $W.$ 

\begin{lemma}\label{lemma1}
For any nonzero $\alpha\in L$ with $|\alpha|^2\neq2$, $E^\alpha \h(-1)=0$.
\end{lemma}  

\begin{proof}
We take an orthonormal basis $\{h_a\}$ so that $h_1\in\C\alpha$. Then $[\w_{a}]h_{b}(-1)\1=\delta_{a,b}h_{b}(-1)\1$.
By \eqref{equation9}, $[E^\alpha]h_{b}(-1)\1=0$ if $b\ne 1.$ In the case
$b=1,$ $\C h_1(-1)$ is a simple module for $A(V_{\Z\a}^+)$
such that $[E^\a]=0$ by Table 2. 
\end{proof}

We next assume that $|\a|^2=2.$ 
Set $\bar{F}^{\alpha}=[E^\alpha]\alpha(-1)\1$.  
Since $[\w_a]\alpha(-1)\1=\alpha(-1)\1$ and $E^\alpha*E^\alpha=4\epsilon(\alpha,\alpha)\w_a,$ we have that $\bar F^\alpha$ is nonzero. 
Using \eqref{equation11} and the fact that $[H_\a]\a(-1)\1=-9\alpha(-1)\1$ 
gives $[H_\a]\bar F^{\alpha}=0$. It is clear that
$\alpha(-1)\1$ and $\bar F^\alpha$ are linear
independent and that $\C \alpha(-1)\1+\C \bar F^{\alpha}$ is closed
under the actions of $[\w_{\alpha}],\,[H_{\alpha}]$ and $[E^{\alpha}]$.  
Since $A(V_{\Z\a}^+)$ is generated by
$[\w_\a],\,[H_\a]$ and $[E^{\alpha}]$ 
we see that $\C \alpha(-1)\1+\C \bar F^{\alpha}$ is an $A(V_{\Z\a}^+)$-module.
By Table 3 we have:
\begin{proposition}\label{prop-3}
For any $\alpha\in L_2$, the two dimensional space $\C \alpha(-1)\1+\C \bar F^{\alpha}$ is an $A(V_{\Z\alpha}^+)$-module isomorphic to $V_{\Z\alpha}^-(0)$. 
\end{proposition}

Now we fix $\a\in L_2$ and prove that $\C\bar F^\a$ is an
irreducible $A(M(1)^+)$-module isomorphic to $M(1,\a)(0).$ 
We continue to fix an orthonormal basis
$\{h_a\}$ of $\h$ such that $h_1\in \C\a$ and  consider the action of 
$A^u,\,A^t$ and $\Lambda_{ab}$ on $\bar F^\a.$ 
By the proof of Lemma \ref{lemma2} we see that 
$[E^t_{aa}]\bar{F}^{\alpha}=[E^\alpha][E^t_{aa}]\alpha(-1)\1=0$ for any $a$.
Thus  $[E^t_{ab}]\bar F^{\alpha}=[E^t_{ab}][E^t_{bb}]\bar{F}^{\alpha}=0$
for any $a,b$ and $A^t\bar F^{\alpha}=0$.

We have already mentioned that
$[H_1]\bar F^{\alpha}=[H_\a]\bar F^{\alpha}=0.$ If $a>1$ then
$[H_a]\bar F^{\alpha}=[E^\a][H_a]\a(-1)\1=0$ by Table
1. Thus by \eqref{equation4}, $[E^u_{aa}]\bar F^{\alpha}=[E^u_{bb}]\bar F^{\alpha}$ for any $a,b.$ 
This shows that $A^u\bar F^{\alpha}=0$ as $[E^u_{ab}]\bar F^{\alpha}=[E^u_{ab}][E^u_{bb}]\bar F^{\alpha}=[E^u_{ab}][E^u_{aa}]\bar F^{\alpha}=0$.

We now deal with $\Lambda_{ab}.$
If $b\neq1$, one has $[E^\alpha]h_b(-1)\1=0$ by \eqref{equation12}. 
Hence $[E^\alpha][E^u_{ba}]h_a(-1)\1=[E^\alpha][E^t_{ba}]h_a(-1)\1=0$.
Lemma \ref{commutelambda} and Table 1 then show that $[\Lambda_{1b}]\bar F^{\alpha}=-[E^\alpha][\Lambda_{1b}]h_a(-1)\1=0$. 
It is clear that $[\Lambda_{ab}]\bar F^{\alpha}=0$ for any $a\neq1,\,b\neq1$.
Consequently, $\C([E^\alpha])\h(-1)=\C\bar{F}^{\alpha}$ is an irreducible $\Free{+}$-module isomorphic to $\Fremo{\alpha}$.

Since $A(V_L^+)$ is generated by $A(M(1)^+)$ and $E^\b$ for $\b\in L,$ 
by Lemma \ref{lemma1}, we see that
\[
W=\h(-1)\bigoplus\left(\sum_{\alpha\in L_2}\C \bar F^{\alpha}\right). 
\]
\begin{lemma}\label{lemmalambda}
Let $f:W\to \charge{-}(0)$ be a linear map defined by $f(h(-1)\1)=h(-1)\1$ and $f(\bar{F}^\alpha)=-2F^\alpha$.
Then $f$ is an $A(\charge{+})$-module isomorphism.
\end{lemma}
\begin{proof}

Clearly, $f$ is well-defined as $\bar F^\a=\bar F^{-\a}$ and
$F^\a=F^{-\a}.$ 
We only have to prove that $f$ is an $A(\charge{+})$-module homomorphism.
We have already proved that $f$ is an $A(\Free{+})$-module homomorphism.
It suffices to show that 
\begin{align}\label{modulehom}
o(E^\beta)f(u)=f([E^\beta]u)
\end{align}
for any $u\in W$ and $\beta\in L$ because $A(\charge{+})$ is generated by $A(\Free{+})$ and $[E^{\beta}]$ for $\beta\in L$. 

Note that $[E^\beta]\h(-1)=0=o(E^\beta)\Free{-}(0)$ if $\beta\notin L_2$. 
So \eqref{modulehom} holds for $u\in \h(-1)$ and $\beta\notin L_{2}$. If $\beta\in L_2$
and $u=h(-1)\1$ such that $(h,\beta)=0,$ then $[E^\b]h(-1)\1=0$ in $W$
and $o(E^\b)h(-1)\1=0$ in $M(1)^-(0).$ Again  \eqref{modulehom} holds
in this case. If $\beta\in L_2$ and $u=\beta(-1)\1,$
\eqref{modulehom} follows from Proposition \ref{prop-3}.
Therefore \eqref{modulehom} holds for $u\in\h(-1)$ and $\beta\in L$.

Now let $u=\bar{F}^{\alpha}$ for some $\alpha\in L_2$. 
Note that $[E^\alpha]\h(-1)=\C \bar{F}^{\alpha}$.
By Proposition \ref{spannalpha1} if $(\beta,\alpha)<0$ then we have 
\begin{align*}
[E^\beta]*[E^{\alpha}]=\sum_{i}[v_{i}]*[E^{\alpha+\beta}]*[w_{i}]
\end{align*}
for some $v_{i},\,w_{i}\in\Free{+}$.
If $\alpha+\beta\notin L_2$ then we have 
$0=f([E^\beta]*\bar{F}^{\alpha})=o(E^\beta)F^\alpha$. 
Hence \eqref{modulehom} hold in this case.
If $\alpha+\beta\in L_2$, $\C[E^{\alpha+\beta}]\h(-1)$ is an $A(\Free{+})$-module isomorphic to $\Fremo{\alpha+\beta}(0)$.
So each $[v_{i}]$ acts as a constant on $\C[E^{\alpha+\beta}]\h(-1)$.
Since \eqref{modulehom} holds for any $u\in\h(-1)$,  
\begin{align*}
f([E^\beta]\bar{F}^{\alpha})&=\sum_{i}f([v_{i}][E^{\alpha+\beta}][w_{i}]\alpha(-1)\1)\\
&=\sum_{i}o(v_{i})o(E^{\alpha+\beta})o(w_{i})f(\alpha(-1)\1)\\
&=\sum_{i}o([v_{i}]*[E^{\alpha+\beta}]*[w_{i}])\alpha(-1)\1\\
&=o([E^\beta])o([E^{\alpha}])\alpha(-1)\1\\
&=-2o([E^\beta])F^{\alpha}\\
&=o([E^\beta])f(\bar{F}^{\alpha}).
\end{align*}
This shows \eqref{modulehom} for $u=\bar{F}^{\alpha}$ and $\beta\in L$ such that $\alpha+\beta\in L_2$.

If $(\alpha,\beta)=0$, we have $E^\beta*E^\alpha=\epsilon(\alpha,\beta)(E^{\alpha+\beta}+E^{\alpha-\beta})$. 
Since $|\alpha\pm\beta|^2=|\alpha|^2+|\beta|^2\geq4$, we see that $o(E^\beta)F^\a=0=[E^\beta]\bar{F}^\alpha$. Thus \eqref{modulehom} holds in this case.
\end{proof}

Thus we get the following proposition:
\begin{proposition}\label{Aunonzro}
Let $W$ be an irreducible $A(\charge{+})$-module such that $A^uW\neq0$.
Then $W\cong\charge{-}(0)$ as $A(\charge{+})$-module.
\end{proposition} 

Finally suppose that $A^tW\neq0$. 
Then we see that $W$ contains an irreducible $A(\Free{+})$-module $\h(-1/2)$ isomorphic to $\Fretw{-}(0)$.
Set $W^0=\{u\in W|A^tu=0\}$. 
If $u\in W^0$, then by Lemma \ref{lemma2}, $A^t[E^\alpha ]u=A^t([I^t][E^\alpha]u)=A^t([E^\alpha][I^t]u)=0$ for any $\a\in L.$ 
That is,
$W^0$ is $[E^\a]$-invariant. Since $A^t$ is a two-sided ideal of
$A(M(1)^+),$ Proposition \ref{spannalpha1} implies that $W^0$ is an $A(\charge{+})$-submodule of $W$. 
Thus $W^0=0$ because $W$ is an irreducible $A(\charge{+})$-module such that $A^tW\neq0$.
Therefore, we have $W=A^tW$. In fact $W$ is a direct sum 
of the unique simple module $\h(-1/2).$

Let $\alpha\in L$ and we take an orthonormal basis $\{h_a\}$ of $\h$ so that $h_1\in\C\alpha$. 
We see that $A^u W=0$ and $[\Lambda_{ab}]W=0$ for any $1\leq a\neq b\leq d$.
Recall the element $[B_\alpha]\in A(\charge{+})$ in \eqref{groupaction}.

\begin{lemma}\label{lemma-Balpha} 
For any $1\leq a,b\leq d$, $[B_{\alpha}]$ and $[E^t_{ab}]$ commute on $W$.
Therefore, 
$[B_\alpha]$ commutes with the action of $A(\Free{+})$.
\end{lemma}

\begin{proof}
It is enough to show the lemma in the case $a=1$ or $b=1$.
Since $A^uW=[\Lambda_{1b}]W=0$, by \eqref{inverse1} and \eqref{inverse2}, we have
\begin{align*}
[E^t_{ab}]=-[S_{ab}(1,1)]-2[S_{ab}(1,2)],\quad [E^t_{ba}]=3[S_{ab}(1,1)]+2[S_{ab}(1,2)]
\end{align*}
on $W$.
Thus \eqref{S1nalpha-equiv} and \eqref{alphaS1n-equiv} gives the following identities; 
\begin{align*}
|\alpha|[E^t_{1b}]*[E^\alpha]=&-4(|\alpha|^2-1)[h_b(-3)F^\alpha]\\
&-(6|\alpha|^2-5)[h_b(-2)F^\alpha]-\left(\frac{3}{2}|\alpha|^2-1\right)[h_b(-1)F^\alpha],\\
|\alpha|[E^t_{b1}]*[E^\alpha]=&4(|\alpha|^2-1)[h_b(-3)F^\alpha]\\
&+(8|\alpha|^2-7)[h_b(-2)F^\alpha]+\left(\frac{9}{2}|\alpha|^2-3\right)[h_b(-1)F^\alpha],\\
|\alpha|[E^\alpha]*[E^t_{1b}]=&-4(|\alpha|^2-1)[h_b(-3)F^\alpha]\\
&-(4|\alpha|^2-5)[h_b(-2)F^\alpha]-\left(\frac{1}{2}|\alpha|^2-1\right)[h_b(-1)F^\alpha],\\
|\alpha|[E^\alpha]*[E^t_{b1}]=&4(|\alpha|^2-1)[h_b(-3)F^\alpha]\\
&+(6|\alpha|^2-7)[h_b(-2)F^\alpha]+\left(\frac{3}{2}|\alpha|^2-3\right)[h_b(-1)F^\alpha].
\end{align*}
These identities implies 
\begin{multline}\label{xxxx}
(2|\alpha|^2-1)[E^t_{1b}]*[E^\alpha]+[E^\alpha]*[E^t_{1b}]\\
=-[E^t_{b1}]*[E^\alpha]
-(2|\alpha|^2-1)[E^\alpha]*[E^t_{b1}].
\end{multline}

Note that $[E^t_{1b}]$ ($[E^t_{b1}]$ resp.\,) is an eigenvector in $A(\charge{+})$ for the left multiplication of $[\w_b]$ of eigenvalue $\frac{1}{16}$ ($\frac{9}{16}$ resp.\,) by \eqref{e-3}. Multiplying 
(\ref{xxxx}) by $[\w_b]$ on the left and using the fact that
$[\w_b]*[E^\a]=[E^\a]*[\w_b]$ we obtain
\begin{multline}\label{xxx}
\frac{1}{16}\!\left((2|\alpha|^2-1)[E^t_{1b}]*[E^\alpha]+[E^\alpha]*[E^t_{1b}]\right)
\\
=-\frac{9}{16}\left([E^t_{b1}]*[E^\alpha]
+(2|\alpha|^2-1)[E^\alpha]*[E^t_{b1}]\right).
\end{multline}
Combining (\ref{xxxx}) and (\ref{xxx}) gives 
\begin{align*}
[E^t_{1b}]*[E^\alpha]&=-\frac{1}{(2|\alpha|^2-1)}[E^\alpha]*[E^t_{1b}]\\
[E^t_{b1}]*[E^\alpha]&=-(2|\alpha|^2-1)[E^\alpha]*[E^t_{b1}].
\end{align*}
It is clear now  that $[E^t_{1b}]*[B_\alpha]=[B_{\alpha}]*[E^t_{1b}]$ and $[E^t_{b1}]*[B_\alpha]=[B_{\alpha}]*[E^t_{b1}]$ by the definition of
$B_\a.$
\end{proof}

Let $\alpha,\beta\in L$ such that $(\alpha,\beta)<0$. 
By Proposition \ref{spannalpha1}, Lemma \ref{lemma-Balpha} and the fact that $W=A^tW$, we can write $[B_\alpha]*[B_\beta]=u*[B_{\alpha+\beta}]$ for some $u\in A^t$ on $W$.
On the other hand $[B_\alpha]*[B_\beta]=u*[B_{\alpha+\beta}]=\epsilon(\alpha,\beta)[B_{\alpha+\beta}]$ on any irreducible $A(\charge{+})$-module $\charge{T_{\chi},-}(0)$. Thus $u=\epsilon(\alpha,\beta)I^t$.
This proves the identity $[B_\alpha]*[B_\beta]=\epsilon(\alpha,\beta)[B_{\alpha+\beta}]$ on $W$ for $\alpha,\beta\in L$ such that $(\alpha,\beta)<0$. 
As in \eqref{zerocase} we can also 
show this identity for $\alpha,\beta\in L$ such that $(\alpha,\beta)=0$.
Then an argument similar to that 
in the proof of Proposition \ref{twistplus} shows that $W$ is an 
$\hat{L}/K$-module on which $\kappa=-1$. 
So $W$ is written as $W=\h(-1/2)\otimes U$ for an $\hat{L}/K$-module $U$. 
By the irreducibility of $W$ we see that $U$ is an
irreducible $\hat{L}/K$-module.   
Hence we have the following proposition:
\begin{proposition}\label{twist-}
Let $W$ be an irreducible $A(\charge{+})$-module such that $A^tW\neq0$. 
Then there exists an irreducible $\hat{L}/K$-module $T_{\chi}$ with central character $\chi$ such that $W\cong\h(-1/2)\otimes T_{\chi}=\charge{T_{\chi},-}(0)$.
\end{proposition}

Combining with Proposition \ref{propositonlambda}, Proposition \ref{twistplus}, Proposition \ref{Aunonzro} and Proposition \ref{twist-}, we see that any irreducible $A(\charge{+})$-module is isomorphic to a top level of known irreducible $\charge{+}$-module. 
Therefore, by Theorem \ref{P3.1} we have the following theorem;

\begin{theorem}\label{classcharge}
Let $L$ be a positive definite even lattice. 
Then any irreducible admissible $\charge{+}$-module is isomorphic to one of irreducible modules $\charge{\pm},\,\charlam{\lambda}\,(\lambda\in L^\circ)$ with $2\lambda\notin L$, $\charlam{\lambda}^{\pm}\,(\lambda\in L^\circ)$ with $2\lambda\in L$ and $\charge{T_{\chi},\pm}$ for any irreducible $\hat{L}/K$-module $T_{\chi}$ with central character $\chi$.
\end{theorem}

\begin{remark}
Although we have achieved the classification of irreducible modules for
$V_L^+,$ we do not determine the structure of $A(V_L^+)$ completely.
By Proposition 5.2 of \cite{Y} and Theorem 5.3 of \cite{ABD}, $V_L^+$ is
$C_2$-cofinite and thus $A(V_L^+)$ is finite dimensional. The orbifold
theory conjectures that $V_L^+$ is a rational vertex operator algebra,
which implies $A(V_L^+)$ is a finite dimensional semisimple associative
algebra. But we cannot prove the semisimplicity of $A(V_L^+)$ in this 
paper.
\end{remark}

\begin{remark} In the case that $L$ is unimodular, $V_L$
has a unique irreducible module $\charge{}$ and a unique $\theta$-twisted module $V_L^{T_{\chi}}.$ So $V_L^+$ has
exactly 4 inequivalent irreducible modules, $V_L^{\pm},$
$(V_L^{T_{\chi}})^{\pm}.$ 
In particular, if $L$ is the Leech lattice, this gives a different proof of
the classification of irreducible $V_L^+$-modules previously obtained 
in \cite{D3}. 
\end{remark}


\begin{thebibliography}{100000}
\bibitem[Y]{Y} 
G. Yamskulna, $C_2$ cofiniteness of vertex operator algebra $V_L^+$ 
when $L$ is a rank one lattice, math.QA/0202056.
\end{thebibliography}
\end{document}